\documentclass[times]{elsarticle}
\usepackage[margin=1.1in]{geometry}

\usepackage{moreverb}
\usepackage[bookmarksopen,bookmarksnumbered,citecolor=red,urlcolor=red]{hyperref}
\usepackage{amsmath}
\usepackage{amsfonts,bm}
\usepackage{graphicx}
\usepackage{amsmath,yhmath}
\usepackage{setspace}
\usepackage{float}
\usepackage{subfigure,epsfig,multirow,bbm,graphicx,graphics, sidecap}
\usepackage{tabularx,ulem}
\usepackage{algorithm} 
\usepackage{algpseudocode}
\usepackage{xcolor}
\usepackage[numbers]{natbib}
\usepackage{algorithmicx}
\usepackage{soul}

\fontencoding{OT1}
\fontfamily{ptm}  
\fontseries{m}    
\fontshape{n}     
\fontsize{11}{11}
\selectfont

\newcommand{\cit}[1]{\citep{#1}}

\newcommand{\M}{\mathbf M}

\newcommand{\x}{ {\mathbf x}}
\newcommand{\y}{ {\mathbf y}}

\newcommand{\R}{\mathbf R}

\newcommand{\A}{\mathbf A}
\newcommand{\B}{\mathbf B}

\newcommand{\mylist}{\begin{list}{$\bullet$}{ \setlength{\topsep}{0cm}
\setlength{\itemsep}{0cm} \setlength{\parsep}{0cm}}}
\newcommand{\mylistend}{\end{list}}

\newcounter{task}
\newcounter{subtask}
\newcounter{subsubtask}

\newcommand{\N}{\mathcal{N}}
\newcommand{\J}{\mathcal{J}}
\renewcommand{\H}{\mathcal{H}}

\newcommand{\V}{\mathbf{V}}
\newcommand{\xa}{\mathbf{x}^\textnormal{a}}
\newcommand{\xb}{\mathbf{x}^\textnormal{b}}

\newcommand{\xtin}{\mathbf{x}_0}
\newcommand{\xk}{\mathbf{x}_k}

\newcommand{\yk}{\mathbf{y}_k}
\newcommand{\nobs}{\textsc{Nobs}}
\newcommand{\nens}{\textsc{Nens}}
\newcommand{\nvar}{\textsc{Nvar}}
\newcommand{\nred}{\textsc{Nred}}

\newcommand{\p}{\mathbf{p}}

\newcommand{\PD}{\mathcal{P}}

\newcommand{\Pa}{\mathcal{P}^{\rm a}}
\newcommand{\tPa}{\widetilde{\mathcal{P}}^{\rm a}}
\newcommand{\dtPa}{{\wideparen{\mathcal{P}}}^{\rm a}}
\newcommand{\Pb}{\mathcal{P}^{\rm b}}
\newcommand{\Bini}{\mathbf{B}_0}
\newcommand{\Aini}{\mathbf{A}_0}
\newcommand{\psdet}{\textnormal{det}^*}
\newcommand{\tr}{\textnormal{Tr}}
\newcommand{\Ss}{\mathbf{S}}

\renewcommand{\Pr}{\textnormal{P}_\textsc{v}}

\renewcommand\bibname{References}

\newtheorem{theorem}{Theorem}[section]
\newtheorem{corollary}{Corollary}[theorem]

\newproof{pf}{Proof}
\newproof{pot}{{\it Proof.}}

\begin{document}
\thispagestyle{empty}
\setcounter{page}{0}

\begin{Huge}
\begin{center}
Computer Science Technical Report CSTR-{1/2016} \\
\today
\end{center}
\end{Huge}
\vfil
\begin{huge}
\begin{center}
{\tt Ahmed Attia, Razvan \c Stef\u anescu, and Adrian Sandu}
\end{center}
\end{huge}

\vfil
\begin{huge}
\begin{it}
\begin{center}
``{\tt The Reduced-Order Hybrid Monte Carlo Sampling Smoother}''
\end{center}
\end{it}
\end{huge}
\vfil

\begin{large}
\begin{center}
Computational Science Laboratory \\
Computer Science Department \\
Virginia Polytechnic Institute and State University \\
Blacksburg, VA 24060 \\
Phone: (540)-231-2193 \\
Fax: (540)-231-6075 \\ 
Email: \url{attia@vt.edu}, \url{sandu@cs.vt.edu} \\
Web: \url{http://csl.cs.vt.edu}
\end{center}
\end{large}

\vspace*{1cm}

\begin{tabular}{ccc}
\includegraphics[width=2.5in]{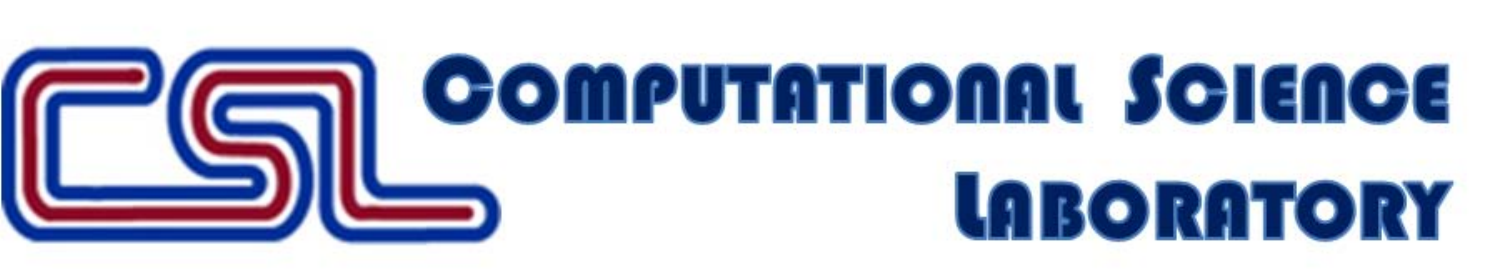}
&\hspace{2.5in}&
\includegraphics[width=2.5in]{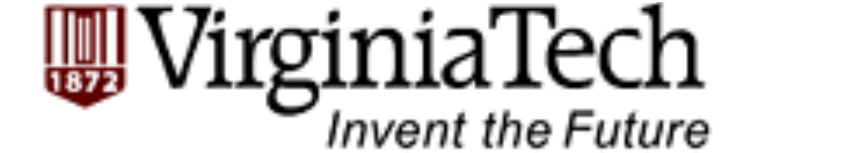} \\
{\bf\em Innovative Computational Solutions} &&\\
\end{tabular}

\newpage

\begin{frontmatter}

\title{The Reduced-Order Hybrid Monte Carlo Sampling Smoother}

\author{A. Attia$^*$, R. \c Stef\u anescu, and A. Sandu }
\address{ Computational Science Laboratory,\\ Department of Computer Science, \\
Virginia Polytechnic Institute and State University, \\
Blacksburg, Virginia, 24060, USA }
\ead{attia@vt.edu,rstefane@vt.edu, sandu@cs.vt.edu}

\cortext[cor1]{Corresponding author}

\date{}
\maketitle

\begin{abstract}
Hybrid Monte-Carlo (HMC) sampling smoother is a fully non-Gaussian four-dimensional data assimilation algorithm that works by directly sampling the posterior distribution formulated in the Bayesian framework. The smoother in its original formulation is computationally expensive due to the intrinsic requirement of running the forward and adjoint models repeatedly. Here we present computationally efficient versions of the HMC sampling smoother based on reduced-order approximations of the underlying model dynamics.
The schemes developed herein are tested numerically using the shallow-water equations model on Cartesian coordinates. The results reveal that the reduced-order versions of the smoother are capable of accurately capturing the posterior probability density, while being significantly faster than the original full order formulation.
\end{abstract}

\end{frontmatter}

\tableofcontents

%
%

\section{Introduction}
\label{Sec:Introduction}
%
Many large-scale prediction problems such as atmospheric forecasting are formulated as initial value problems. The uncertainty of the associated model initial conditions can be decreased by combining imperfect forecasts produced by propagating the model dynamics with real measurements collected at discrete times over an assimilation window. The model state and the observations are both uncertain and can be modeled as random variables. The probability distribution describing the knowledge about the initial state system, before incorporating observations information, is known as the prior distribution. The likelihood function accounts for the discrepancies between the measurements and model-predicted observations. Data assimilation applies Bayes' theorem to obtain a posterior probability distribution, named the analysis, that characterizes the knowledge about the system state given the observed measurements. In practice, due to the state space high-dimensionality of many realistic models, it is impossible to exactly describe the posterior distribution and several assumptions and approximations are necessary. Widely accepted assumptions are that the background and the observation errors are characterized by Gaussian distributions, with no correlations between observation errors at different time instances.

Two families of methodologies are generally followed in order to generate accurate estimates of the true system state.   Ensemble-based statistical methods seek to approximate the posterior probability density function (PDF) based on an ensemble of model states, while the variational approaches estimate the true state of the system by searching for the state that maximizes the posterior PDF. Among the variational techniques, the 4D-Var method achieves this goal by searching for a local minimum of an objective function corresponding to the negative logarithm of the posterior distribution.  4D-Var finds the maximum posterior (MAP) estimate of the true state and does not directly estimate the uncertainty associated with the analysis state.  For scenarios including nonlinear observation operators and state models, the analysis distribution is not Gaussian, and the 4D-Var algorithm may be trapped in a local minimum of the cost function leading to incorrect conclusions.

Accurate solution of the non-Gaussian data assimilation problems requires accounting for all regions of high probability in the posterior distribution. The recently developed hybrid Monte-Carlo (HMC) sampling smoother~\cite{attia2014HMCsampling,Attia_HMCSmoother_TR} is a four dimensional data assimilation scheme designed to solve the non-Gaussian smoothing problem by sampling from the posterior distribution. It relies on an accelerated Markov chain Monte-Carlo (MCMC) methodology where a Hamiltonian system is used to formulate proposal densities. Two issues are important when using HMC sampling strategy. First, it requires the formulation of the posterior negative-log function gradient.
Secondly, the involved Hamiltonian system is propagated using a symplectic integrator whose parameters must be carefully tuned in order to achieve good performance.

While producing consistent description of the updated system uncertainty (e.g., the analysis error covariance matrix), the original formulation of the HMC sampling smoother is computationally expensive when compared to the 4D-Var approach. This is due to the large number of gradient evaluations required, which translates into many forward and adjoint models runs. In the case of large scale problems the computational cost becomes prohibitive. In this present study we propose a practical solution by approximating the gradient using information obtained from lower-dimensional subspaces via model reduction.

Reduced order modeling refers to the development of low-dimensional systems that represent the important  characteristics of a high-dimensional or infinite dimensional dynamical system. Typically this is achieved by projecting model dynamics onto a lower dimensional spaces. Construction of low relevant manifolds can be achieved using the reduced basis method \citep{BMN2004,grepl2005posteriori,patera2007reduced,rozza2008reduced,Dihlmann_2013,Lieberman_et_al_2010}, dynamic mode decomposition \citep{Rowley2009,Schmid2010,Tissot_et_al_2014,Bistrian_Navon_2014} and Proper Orthogonal Decomposition (POD) \citep{karhunen1946zss,loeve1955pt,hotelling1939acs,lorenz1956eof}. The latest is the most prevalent basis selection method for nonlinear problems.  Data analysis using POD and method of snapshots \citep{Sir87a, Sir87b, Sir87c} is conducted to extract basis functions, from experimental data or detailed simulations of high-dimensional systems, for subsequent use in Galerkin projections that yield low dimensional dynamical models. By coupling POD with empirical interpolation method
(EIM) \citep{BMN2004}, discrete variant DEIM \citep{Cha2008,ChaSor2012,ChaSor2010,Drmac_Gugercin2015} or best points interpolation method \citep{NPP2008}, one can obtain fast approximations of reduced order nonlinear terms independent of the dimension of high-fidelity space.  Other such approaches include missing point estimation \citep{Astrid_2008} and Gauss-Newton with approximated tensors \citep{Carlberg2_2011,Carlberg_2012} relying upon the gappy POD technique \citep{Everson_1995}.

The present manuscript develops practical versions of the HMC sampling smoother using approximate dynamical information obtained via POD/DEIM reduced order models \citep{Stefanescu2013,Stefanescu_etal_forwardPOD_2014,dimitriu2014application}. Two reduced order variants are proposed differentiated by the choice of the sampling space used to generate the proposals. In the first case we are sampling in a reduced order space while the second approach samples directly from the high fidelity space. For both scenarios the negative logarithm of the posterior distribution reassembles one of the flavors of reduced order 4D-Var objective functions \citep{fang2009pod,cao2007reduced,vermeulen2006model}, where the prior term is either estimated in the reduced or full space respectively. While reduced order models are employed to estimate the posterior distribution negative logarithm and it's gradient, the associated Hamiltonian system generates proposals in a reduced or high-fidelity space depending on the nature of the sampling space.

The choice of reduced order manifolds is crucial for the accuracy of the reduced samplers. The smoothers may suffer from the fact that the basis elements are computed from a reference trajectory containing features which are quite different from those of the current proposal trajectory. To overcome the problem of unmodelled dynamics in the POD-basis we propose to update the basis from time to time according to the current state similarly as in the case of adaptive reduced order optimization \cite{Afanasiev_Hinze_2001,Ravindran_2002,Kunisch_2008}. Inspired by the recent advances in reduced order data assimilation field where it was shown that accurate reduced order Karush-Kuhn-Tucker conditions with respect to their full order counterparts highly increase the accuracy of the obtained analysis \cite{cstefuanescu2015pod}, we chose to update the reduced order manifolds based on high-fidelity forward and adjoint trajectories corresponding to current proposal as well as the associated gradient of the negative logarithm of the high-fidelity posterior distribution. The basis is refreshed once after several HMC smoothers iterations. The numerical results using swallow-water equations model on cartesian coordinates reveal
that the reduced-order versions of the smoother are accurately capturing the posterior probability
density, while being significantly faster than the original full order formulation.

The paper is organized as follows. Section~\ref{Sec:DA} reviews the four-dimensional data assimilation problem and the original HMC smoother. Section~\ref{Sec:Reduced_Order_DA} reviews reduced-order modeling, and introduces the reduced order 4D-Var data assimilation framework. Section~\ref{Ssec:Reduced_HMC_Smoother} presents a rigorous description of the proposed versions of the HMC smoother. In section~\ref{Sec:Reduced_Order_Sampling} several theoretical properties of the distributions sampled with reduced order models are derived. Numerical results are presented in Section~\ref{Sec:Numerical_Results} while conclusions are drawn in Section~\ref{Sec:Conclusion}.

\section{Data Assimilation}
\label{Sec:DA}

One can solve a data assimilation (DA) problem by describing the posterior distribution in the Bayesian formalism given the prior and the likelihood function. In the four-dimensional DA (4DDA) context, the main goal is to describe the posterior distribution of system state at the initial time of a specific assimilation window. The assimilation window is defined based on the availability of observations at specific discrete time instances. The knowledge about the system at the initial time $t_0$ defines the distribution (the prior) $\Pb(\xtin)$ of the state $\xtin \in \mathbb{R}^{\nvar}$ before absorbing knowledge presented by the observations captured over the assimilation window.
Based on a set of observations $\{\y_k=\y[t_k] \in \mathbb{R}^{m}\}_{k=0,1,\ldots,\nobs}$, at the discrete time points $\{t_k\}_{k=0,1,\ldots,\nobs}$ in the interval $[t_0\,,t_F]$, the sampling distribution (likelihood function) is defined as $\PD(\y_0,\y_1,\ldots,\y_{\nobs} | \xtin)$, and the posterior distribution describing the updated information about the model state at the initial time, given the measurements, takes the general form
\begin{equation}
\label{eqn:Bayes_Rule}
\PD(\xtin | \y_0,\y_1,\ldots,\y_{\nobs}) = \frac{\Pb(\xtin)\, \PD(\y_0,\y_1,\ldots,\y_{\nobs} | \xtin))}{\PD(\y_0,\y_1,\ldots,\y_{\nobs})}\,.
\end{equation}
Fully and accurately describing this general posterior probability density function in DA literature is an intractable problem, and usually simplifying assumptions are generally made. As we mentioned in Section~\ref{Sec:Introduction}, a frequent supposition is that the background and observation errors are characterized by Gaussian distributions. If the observations are assumed to be independent from the model states, and the associated error characteristics of these observations are not correlated in time, the posterior distribution takes the form
\begin{subequations}
\label{eqn:Smoothing_Posterior_Kernel}
\begin{align}
\Pa(\xtin) &= \PD(\xtin | \y_0,\y_1,\ldots,\y_{\nobs}) \propto \exp{ \Bigl( -\J(\xtin) \Bigr)}\,, \label{eqn:Full_kernel_distribution}\\
\J:\mathbb{R}^{\nvar} \to \mathbb{R},\quad \J(\xtin) &= \frac{1}{2} \lVert \xtin-\xb_0 \rVert ^2 _{ \mathbf{B}_0^{-1} }
+ \frac{1}{2}\sum_{k=0}^{\nobs}{ \lVert \y_k - \mathcal{H}_k(\x_k) \rVert^2 _{ \mathbf{R}_k^{-1} } } \,, \label{eqn:4DVar_Cost_Functional}
\end{align}
\end{subequations}
where $\xb_0 =\xb[t_0]$ is a background/forecast state, the matrices $ \mathbf{B}_0$ and $\mathbf{R}_k,~{k=0,1,\ldots,\nobs}$, are the covariance matrices associated with the background and
measurement errors respectively. The observation operator $\mathcal{H}_k : \mathbb{R}^{\nvar} \to \mathbb{R}^{m},$ at time instance $t_k$, maps a given state $\x_k \in \mathbb{R}^{\nvar}$ to the
observation space. The dimension of the observation space is usually much smaller than the size of the state space, that is $m \ll \nvar$. The associated norms over the state and observation
spaces are defined as:
\begin{equation}
\label{eqn:norm_def}
\| \mathbf{a}\,-\mathbf{b}\|^2_\mathbf{C} = (\mathbf{a} - \mathbf{b})^T\mathbf{C}(\mathbf{a} - \mathbf{b}),
\end{equation}
where $\mathbf{a},~\mathbf{b}$ belong to either $\mathbb{R}^{\nvar}$ or $\mathbb{R}^m$ and $\mathbf{C}$ is a matrix of $\mathbb{R}^{\nvar \times \nvar}$ or $\mathbb{R}^{m \times m}$ dimensions.

The model state $\x_k=\x[t_k]$, is obtained from the initial state $\xtin$, by propagating the model dynamics to time point $t_k$, i.e.
\begin{equation}
\label{eqn:Forward_Model}
\x_k = \mathcal{M}_{t_0\rightarrow t_k}(\xtin)= \mathcal{M}_{0,k}(\xtin)\,,
\end{equation}
where, $\mathcal{M}_{t_0\rightarrow t_k} : \mathbb{R}^{\nvar} \to \mathbb{R}^{\nvar},~k=1,\ldots,\nobs$ represents the discretized mathematical model reflecting the state dynamics. The function $\J$ given in~\eqref{eqn:Smoothing_Posterior_Kernel} is quadratic, and consequently the distribution~\eqref{eqn:Smoothing_Posterior_Kernel} is Gaussian distribution, only if the states $\xk$ are linearly related to the system initial condition $\xtin$, and the observations $\yk$ are linearly related to the model states $\xk$. In virtually all practical settings, the time dependent models are complicated resulting in a nonlinear bond between
states at different time instances. More difficulty is added since for example, in the field of atmospheric sciences, highly-nonlinear observation operators are often constructed to relate new types of measurements to state variables.

\subsection{4D-Var data assimilation}
\label{subsec:Full_4DVAR_DA}
The strongly-constrained 4D-Var DA scheme is an optimization algorithm that searches for the maximum aposteriori estimate (MAP) by seeking a \textit{local} minimizer of the cost function~\eqref{eqn:4DVar_Cost_Functional}, constrained by the model equation~\eqref{eqn:Forward_Model}. Precisely, the constrained optimization problem is defined as
\begin{equation}
\label{eqn:4DVAR_Optimization}
\begin{aligned}
\min_{\xtin}{
  \mathcal{J}(\xtin) } &= \frac{1}{2} \lVert \xtin-\xb_0 \rVert ^2 _{ \mathbf{B}_0^{-1} }
+ \frac{1}{2}\sum_{k=0}^{\nobs}{ \lVert \y_k - \mathcal{H}_k(\x_k) \rVert^2 _{ \mathbf{R}_k^{-1} } }   \,, \\
\x_k &= \mathcal{M}_{t_0\rightarrow t_k}(\xtin)\,~k=1,\ldots,\nobs.
\end{aligned}
\end{equation}
Gradient-based schemes are generally employed to find the corresponding initial conditions which in the control theory language are known as control variables. Using the Lagrange multiplier technique the constrained optimization problem \eqref{eqn:4DVAR_Optimization} is  redesigned as an unconstrained minimization problem with the associated first order necessary optimality conditions:
\begin{subequations}
\label{eqn:cici}
\begin{eqnarray}
\label{eqn:KKT_Full_forward_model}
 \textnormal{{\it Forward model}:  }&&
~~\x_k = \mathcal{M}_{t_0\rightarrow t_k}(\xtin),~k=1,\ldots,\nobs, \quad \quad \quad \quad \quad \\
\label{eqn:Full_adjoint_model}
 \textnormal{{\it Adjoint model}:  }  && \boldsymbol{\lambdaup}_{\nobs} = {\mathbf H}_{\nobs}^T \mathbf{R}_{\nobs}^{-1} \big({\bf y}_{\nobs} - \H_{\nobs}({\mathbf x}_{\nobs})\big), \quad \quad \quad \quad \quad \\
&& \boldsymbol{\lambdaup}_{k} = \mathbf{M}^T_{t_0\rightarrow t_{k+1}}{\boldsymbol{\lambdaup}}_{k+1} + \mathbf{H}_k^T \mathbf{R}_{k}^{-1} \big({\bf y}_k - \H({\mathbf x}_k)\big),\quad k=\nobs-1,..,0, \\
\label{eqn:Full_function_gradient}
\textnormal{{\it Cost function gradient:} } && \nabla_{{\bf x}_0}\mathcal{J}(\xtin) = -{\mathbf B}_0^{-1}\big({\mathbf x}^{\rm b}_0-{\mathbf x}_0\big) - {\boldsymbol{\lambdaup}}_0 = 0. \quad \quad \quad \quad \quad
\end{eqnarray}
\end{subequations}

The Jacobians of the model and observation operators are denoted by $\mathbf{M}_{t_0\rightarrow t_{k+1}},~k=0,1,\ldots,\nobs-1$ and $\mathbf{H}_k,~k=0,1,\ldots,\nobs$ while the adjoint solution $\boldsymbol{\lambdaup}_k \in \mathbb{R}^{\nvar},~k=0,1,\ldots,\nobs$, provides an efficient way to compute the gradient~\eqref{eqn:Full_function_gradient}.
The nonlinear optimization procedure is computationally expensive  and either low-resolution models (incremental 4D-Var~\cite{Stefanescu_lightning_2013}), or alternatively reduced-order models are used~\cite{cstefuanescu2015pod} to alleviate this drawback.
It is well-known that 4D-Var does not inherently provide a measure of the uncertainty about the updated state (e.g., analysis error covariance matrix) and usually hybrid methods are
considered to account for this type of information. This approach results in some inconsistency between the analysis state and the analysis error covariance matrix especially when they are obtained using different algorithms.

%
\subsection{Smoothing by sampling and the HMC sampling smoother}
\label{Subsec:HMC_Sampling_Smoother}
Monte-Carlo smoothing refers to the process of representing/approximating the posterior distribution~\eqref{eqn:Smoothing_Posterior_Kernel} using an ensemble of model states sampled from that posterior. Ensemble Kalman smoother (EnKS)~\cite{Evensen_2000} is an extension of the well-known ensemble Kalman Filter~\cite{Kalman_1960} to the case where observation are assimilated simultaneously. EnKS produces a minimum-variance unbiased estimate (MVUE) of the system state by estimating the expectation of the posterior $\mathtt{E}_{\Pa}{[\xtin]}$ using the mean of an ensemble of states. The strict Gaussianity and linearity assumptions imposed by EnKS, usually result in poor performance of the smoother.

Pure sampling of the posterior distribution~\eqref{eqn:Bayes_Rule} using a Markov Chain Monte-Carlo (MCMC)~\cite{metropolis1953equation,neal2011mcmc} technique is known in theory to provide more accurate estimates without strictly imposing linearity or Gaussianity constraints. MCMC is a family of Monte-Carlo schemes tailored to sample a given distribution (up to a proportionality constant), by constructing a Markov chain whose stationary distribution is set as the target distribution. By design, an MCMC sampler is guaranteed to converge to its stationarity. However, the choice of the proposal density, the convergence rate, the acceptance rate, and the correlation level among sampled points are the main building blocks of MCMC responsible for its performance and efficiency. Practical application of MCMC requires developing accelerated chains those can attain stationarity fast, and then explore the state space efficiently in very few steps. One of the MCMC samplers mainly designed for complicated PDFs and large dimensional spaces is the Hamiltonian/Hybrid Monte-Carlo sampler (HMC). HMC was firstly presented in~\cite{duane1987hybrid} as an accelerated MCMC sampling algorithm. The sampler mainly uses information about the geometry of the posterior to guide its steps in order to avoid random walk behaviour and visit more frequently regions with high probability with the capability of jumping between separated modes of the target PDF.

\paragraph{HMC sampling}
%
As all other MCMC samplers, HMC samples from a PDF
\begin{equation}
\label{eqn:HMC_PDF_kernel}
\pi{(\x)}\propto \exp{(-\J^N(\x))};\, \x \in \mathbf{R}^{\nvar}\,,
\end{equation}
where $\exp{(-\J^N(\x))}$ is the shape function of the distribution, and $\J^N:\mathbb{R}^{\nvar} \to \mathbb{R}$ is the PDF negative-log.
The power of MCMC, and consequently HMC, is that only the shape function, or alternatively the negative-log, is needed, while the scaling factor is not strictly required as in the case of standard application of Bayes' theorem.
HMC works by viewing $\x$ as a position variable in an extended phase space consisting of points $(\p,\x)\, \in \mathbf{R}^{2\nvar}$, where $\p \in \mathbb{R}^{\nvar}$ is an auxiliary momentum variable.
The Hamiltonian dynamics is modeled by the set of ordinary differential equations (ODEs):
\begin{equation}
\label{eqn:Hamiltonian_equations}
\begin{aligned}
 \frac{d\x}{dt} &= \nabla_\p\, H\,, \\
 \frac{d\p}{dt} &= - \nabla_\x\, H\,,
\end{aligned}
\end{equation}
where $H=H(\p,\x)$ is the constant total energy function of the system, known as the the Hamiltonian function or simply the Hamiltonian. A standard formulation of the Hamiltonian in the context of HMC is:
\begin{equation}
\label{eqn:Hamiltonian_function}
H(\p,\x)  =  \underbrace{\frac{1}{2} \, \p^T \mathbf{M}^{-1} \p}_{\text{ kinetic energy}} + \underbrace{\mathcal{J}^N(\x)}_{\text{ potential energy}}\,,
\end{equation}
where $\M \in \mathbb{R}^{\nvar \times \nvar}$ is a positive definite matrix known as the mass matrix.
This particular formulation leads to a canonical distribution of the joint state $(\p,\, \x)$ proportional to:
\begin{equation}
\label{eqn:Canonical_PDF_Kernel}
 \exp{(-H(\p,\, \x))} = \exp{\left( - \frac{1}{2} \, \p^T \mathbf{M}^{-1} \p \right)}\,\, \pi(\x)\,.
\end{equation}

The exact flow $\Phi_T:\mathbb{R}^{2\nvar} \rightarrow \mathbb{R}^{2\nvar}; \quad \Phi_T\bigl(\p[0],\x[0]\bigr)=\bigl(\p[T],\x[T]\bigr)\,$ describes the time evolution of the Hamiltonian system~\eqref{eqn:Hamiltonian_equations} and is practically approximated using a numerical integrator that is symplectic as
$\phi_T:\mathbb{R}^{2\nvar} \rightarrow \mathbb{R}^{2\nvar}; \quad \phi_T\bigl(\p[0],\x[0]\bigr)\approx \bigl(\p[T],\x[T]\bigr)\,$.
The use of a symplectic numerical integrator to approximate the exact Hamiltonian flow results in changes of total energy.
Traditional wisdom recommends splitting the pseudo-time step $T$ into $m$ smaller steps of size $h$ in order to simulate the Hamiltonian trajectory between points $\bigl(\p[0],\x[0]\bigr),\, \bigl(\p[T],\x[T]\bigr)$ more accurately.
The symplectic integrator of choice is position (or velocity) Verlet. New higher order symplectic integrators were proposed recently and tested in the context of filtering for data assimilation~\cite{sanz2014markov,Attia_HMCFilter_TR}.
One step of size $h$ of the position Verlet \cite{sanz1994numerical,sanz2014markov} integrator is describes as follows
\begin{subequations}
\label{eqn:Position_Verlet}
\begin{align}
\x[{h/2}]  &=  \x[0] + \frac{h}{2}\, \mathbf{M}^{-1}\, \p[0] \,,   \\
 \p[{h}]   &=  \p[0] - h \,\nabla_\x \mathcal{J}(\x[{h/2}])  \,,   \label{eqn:Position_Verlet_Middle}\\
 \x[{h}]   &=  \x[{h/2}] + \frac{h}{2}\, \mathbf{M}^{-1} \,\p[{h}].
\end{align}
\end{subequations}
While, the mass matrix is a user-defined parameter, it can be designed to enhance the performance the sampler~\cite{Attia_HMCFilter_TR}.
The step parameters of the symplectic integrator $m,\, h$ can be empirically chosen by monitoring the acceptance rate in a preprocessing step.
Specifically, the parameters of the Hamiltonian trajectory can be empirically adjusted such as to achieve a specific rejection rate.
Generally speaking, the step size should be chosen to achieve a rejection rate between $25\%$ and $30\%$, and the number of steps should generally be large~\cite{neal2011mcmc}.

Adaptive versions of HMC have been also proposed with the capability of adjusting it's step parameters. No-U-Turn sampler (NUTS)~\cit{hoffman2014NUTS} is a version of HMC capable of automatically tuning its parameters to prohibit the sampler from retracing its steps along the constructed Hamiltonian  trajectory.
Another HMC sampler that tunes its parameters automatically using third-order derivative information is the Riemann manifold HMC (RMHMC)~\cite{girolami2011riemann}.

The intuition behind HMC sampler is to build a Markov chain whose stationary distribution is defined by the canonical PDF~\eqref{eqn:Canonical_PDF_Kernel}. In each step of the chain,
a random momentum $\p$ is drawn from a Gaussian distribution $\mathcal{N}(0,\M)$, and the Hamiltonian dynamics \eqref{eqn:Hamiltonian_equations} at the final pseudo-time interval
proposes a new point that is either accepted or rejected using a Metropolis-Hastings criterion. The two variables $\p,\, \text{and } \x$ are independent, so discarding the momentum
generated at each step will leave us with sample points $\x$ generated from our target distribution.

In a previous work~\cite{Attia_HMCSmoother_TR} we proposed using HMC as a pure sampling smoother to solve the nonlinear 4DDA smoothing problem. The method samples from the posterior
distribution of the model state at the initial time on an assimilation window on which a set of observations are given at discrete times. Following general assumptions where the prior is
Gaussian and the observation errors are normally distributed, the target distribution defined in \eqref{eqn:HMC_PDF_kernel} is identical with the posterior
distribution associated with the smoother problem \eqref{eqn:Smoothing_Posterior_Kernel}. Consequently the PDF negative-log $\J^N$ in \eqref{eqn:HMC_PDF_kernel} resembles the 4D-Var cost
function $\J$ defined in \eqref{eqn:4DVAR_Optimization} and the gradient of the potential energy required by the symplectic integrator is the gradient of the 4D-Var cost functional
\eqref{eqn:Full_function_gradient}.
The main hindrance stems from the requirement of HMC to evaluate the gradient of the potential energy (target PDF negative-log) at least as many times as the symplectic integrator is involved,
which is an expensive process.
Despite the associated computational overhead, the numerical results presented in~\cite{Attia_HMCSmoother_TR} show the potential of using HMC smoother to sample multi-modal, high-dimensional
posterior distributions formulated in the smoothing problem.

\section{Four-Dimensional Variational Data Assimilation with Reduced-Order Models}
\label{Sec:Reduced_Order_DA}
Optimization problems such as the one described in \eqref{eqn:4DVAR_Optimization} for nonlinear partial differential equations often demand very large computational resources, so that the need for developing fast novel approaches emerges. Recently the reduced order approach applied to optimal control problems for partial differential equations has received increasing attention. The main idea is to project the dynamical system onto subspaces consisting of basis elements that represent the characteristics of the expected solution. These low order models serve as surrogates for the dynamical system in the optimization process and the resulting approximate optimization problems can be solved efficiently.

\subsection{Reduced order modeling}
\label{Subsect:ROM}
Reduced order modeling refers to the development of low-dimensional models that represent desired characteristics of a high-dimensional or infinite dimensional dynamical system.
Typically, models are constructed by projection of the high-order, high-fidelity model onto a suitably chosen low-dimensional reduced-basis \cite{Antoulas_2005_book}. Most reduced-bases for nonlinear problems are constructed from a collection of simulations (methods of snapshots \cite{Sir87a,Sir87b,Sir87c})

The most popular nonlinear model reduction technique is Proper Orthogonal Decomposition (POD) and it usually involves a Galerkin projection with basis $V \in \mathbb{R}^{\nvar \times \nred}$
obtained as the output of Algorithm \ref{alg:svd_basis}. Here $\nred$ is the dimensional of the reduced-order state space spanned e.g., by the POD basis.

\begin{algorithm}
 \begin{algorithmic}[1]
 \State Solve for the state variable solutions $\x_k,~k=1,..,\nobs$ of \eqref{eqn:Forward_Model}. One can make use of more snapshots to construct the basis thus for example to consider a number of time steps larger than $\nobs$.
 \State Compute the singular value decomposition (SVD) for the state variable snapshots matrix $[{\bf x}_{0}\,{\bf x}_{1}~~ \ldots ~~{\bf x}_{\nobs}] = \bar {\mathbf{V}} \Sigma {\bar W}^T,$
with the singular vectors matrix $\bar {\mathbf{V}} =[{\bf v}_i]_{i=1,..,\nvar}$.
 \State Using the singular-values $\lambda_1\geq \lambda_2\geq ...~\lambda_n\geq 0$ stored in the diagonal matrix $\Sigma$,
define $I(p)=( {\sum_{i=1}^p \lambda_i)/(\sum_{i=1}^{\nvar} \lambda_i})$.
\State Choose $\nred$, the dimension of the POD basis, such that $ \nred=\min_p \{I(p):I(p)\geq \gamma\}$ where $0 \leq \gamma \leq 1$ is the percentage of total information captured by the reduced space $\mathcal{X}^{\nred}=\textnormal{range}(\mathbf{V})$, usually $\gamma=0.99$.
 \end{algorithmic}
 \caption{POD basis construction}
 \label{alg:svd_basis}
\end{algorithm}

Assuming a POD expansion $\x_k \approx \V \tilde{\x}_k,~\tilde{\x}_k \in \mathbb{R}^{\nred},~k=0,..,{\nobs}$ (for simplicity we neglected the centering trajectory, shift mode or mean field
correction~\cite{NAMTT03}) and making use of the basis orthogonality the associated POD-Galerkin model of~\eqref{eqn:Forward_Model} is obtained as
\begin{equation}
\label{eqn::POD_Galerkin_model}
\tilde{\x}_{k+1} = \mathbf{V}^T\mathcal{M}_{t_{k}\to t_{k+1}}(\mathbf{V}\tilde{\x}_k),~k=0,..,\textsc{N}_{\nobs}.
\end{equation}
The efficiency of POD - Galerkin technique is limited to the linear or bilinear terms \cite{stefanescu2014comparison} and strategies such as Empirical Interpolation Method (EIM)
\cite{MBarrault_YMaday_NDNguyen_ATPatera_2004a}, Discrete Empirical Interpolation Method (DEIM)~\cite{ChaSor2010,Stefanescu2013} and tensorial POD \cite{stefanescu2014comparison}are
usually employed to alleviate this deficiency.

\subsection{Reduced order 4D-Var data assimilation}
\label{Subsec:RO_DA}
%

The ``adjoint of reduced plus reduced of adjoint'' approach (ARRA) leads to the construction of consistent feasible reduced first order optimality conditions \cite{cstefuanescu2015pod} and this framework is employed to build the reduced POD manifolds for the reduced HMC samplers. In the case of Galerkin projection the POD reduced space is constructed based on sampling of both full forward and adjoint trajectories as well as the gradient of the cost function background term.

The reduced data assimilation problem minimizes the following reduced order cost function $\J^\textsc{pod} : \mathbb{R}^k \to \mathbb{R}$
\begin{subequations}
\label{eqn:reduced_optimization}
\begin{equation}\label{eqn:reduced_cost_func}
\J^\textsc{pod}(\tilde{\x}_0)= \frac{1}{2} \lVert \V \tilde{\x}_0-\xb_0 \rVert ^2 _{ \Bini^{-1} }
+ \frac{1}{2}\sum_{k=0}^{\nobs}{ \lVert \y_k
  - \H_k(\V\tilde{\x}_k) \rVert^2 _{ \mathbf{R}_k^{-1} } } \,,
\end{equation}
subject to the constraints posed by the ROM projected nonlinear forward model dynamics~\eqref{eqn::POD_Galerkin_model}
\begin{equation}\label{eqn:reduced_forward_AR}
\tilde{\x}_{k+1} = \widetilde{\mathcal{M}}_{t_{k} \to t_{k+1}}(\tilde{\x}_k),\,
\widetilde{\mathcal{M}}_{t_{k} \to t_{k+1}}(\tilde{\x}_k) = \V^T \mathcal{M}_{t_{k} \to t_{k+1}}(\V {\bf \tilde x}_k),\, k=0,1,\ldots,\textsc{N}_{\nobs}\,.
\end{equation}
\end{subequations}
An observation operator that maps directly from the reduced model space to observations space may be introduced, however for this study the operator requires the projected states as
input.

The associated reduced Karush-Kuhn-Tucker conditions \cite{cstefuanescu2015pod} are:
\begin{subequations}
\label{eqn:KKT_AR}
\begin{eqnarray}
\label{eqn:reduced_KKT_forward_AR}
&& \textnormal{{\it ARRA reduced forward model}:  }\\
\nonumber
&& \qquad \tilde{\x}_{k+1} = \widetilde{\mathcal{M}}_{t_{k} \to t_{k+1}}\left(\tilde{\x}_k\right),\quad k=0,..,\textsc{N}_{\nobs}\,, \\
\label{eqn:adjoint_reduced_model}
&& \textnormal{{\it ARRA reduced adjoint model}:  }  \\
\nonumber
&& \qquad \tilde{\boldsymbol{\lambdaup}}_{\nobs} = \V^T\widehat{\mathbf{H}}_{\nobs}^T \mathbf{R}_{\nobs}^{-1} \big({\bf y}_{\nobs} - \H_{\nobs}(\V\tilde{\x}_{\nobs})\big), \quad \quad \quad \quad \quad \quad \quad \\
\nonumber
&& \qquad \tilde{\boldsymbol{\lambdaup}}_k = \V^T\widehat{\mathbf{M}}^T_{t_{k} \to t_{k+1}}\V\tilde{\boldsymbol{\lambdaup}}_{k+1} + \V^T\widehat{\mathbf{H}}_k^T \mathbf{R}_{k}^{-1} \big({\mathbf y}_k - \H_k(\V\tilde{\x}_k)\big),\quad k=\nobs-1,..,0, \\
\label{eqn:reduced_cost_function_gradient}
&& \textnormal{{\it ARRA cost function gradient} }: \\
\nonumber
&& \qquad \nabla_{\tilde{\x}_0}\mathcal{J}^\textsc{pod} = -\V^T {\bf B}_0^{-1}\big({\mathbf x}^{\rm b}_0-\V\tilde{\x}_0\big) - \tilde{\boldsymbol{\lambdaup}}_0 = 0.
\end{eqnarray}
\end{subequations}
The operators $\widehat{\mathbf{H}}_k,~k=0,..,\nobs$ and $\widehat{\mathbf{M}}_{t_{k} \to t_{k+1}},~k=0,..,\nobs$ are the Jacobians of the high-fidelity observation operator $\mathcal{H}_k$ and model $\mathcal{M}_{t_{k} \to t_{k+1}}$ evaluated at $V \tilde{\x}_k$.

\section{Reduced-Order HMC Sampling Smoothers}
\label{Ssec:Reduced_HMC_Smoother}
%
One of the key features of HMC sampler is the clever exploration of the state space with guidance based on the distribution geometry. For this the HMC sampling smoother requires not only forward model propagation to evaluate the likelihood term, but also the evaluation of the gradient of the negative-log of the posterior PDF using the adjoint model. This gradient can be approximated using information from the reduced space as obtained from the reduced order model \eqref{eqn:reduced_forward_AR}.

The HMC algorithm requires that both the momentum $\p$ and the target state $\x$ are vectors of the same dimension. There are two ways to achieve this while using reduced order information:
\begin{enumerate}
\item[i) ] sample the reduced-order subspace only, i.e., collect samples from~\eqref{eqn:Fully_Projected_PDF}, or
\item[ii)] sample the full space, i.e., collect samples from~\eqref{eqn:Full_PDF_with_Approximate_Gradient} but use an approximate gradient of the posterior negative-log likelihood function obtained in the reduced space.
\end{enumerate}
We next discuss each of these options in detail.

\subsection{Sampling in the reduced-order space}
In this approach the model states are fully projected in the reduced space, $\tilde{\x}\in \mathbb{R}^{\nred}$. The target posterior distribution,
and the potential energy are given by:
\begin{subequations}
\label{eqn:Reduced_HMC_Potential_Energy_and_Gradient}
\begin{equation}
\label{eqn:Fully_Projected_PDF}
\begin{aligned}
  \pi{(\tilde{\x}_0)} &\propto \exp{\big(-\widetilde{\J}(\tilde{\x}_0)\big)}\,, \\
  \widetilde{\J}(\tilde{\x}_0) &= \frac{1}{2} \big\lVert \tilde{\x}_0 - \V^T\xb_0 \big\rVert^2 _{\left(\V^T \B_0 \V \right)^{-1} }
  + \frac{1}{2}\sum_{k=0}^{\nobs}{ \big\lVert \y_k - \H_k(\V\tilde{\x}_k) \big\rVert^2 _{ \mathbf{R}_k^{-1} } } \,, \\
  \tilde{\x}_k &= \widetilde{\mathcal{M}}_{k-1,\,k}(\tilde{\x}_{k-1}),\,\quad k=1,0,\ldots,\nobs\,; \quad \tilde{\x}_0 = \V^T \xtin \,,
\end{aligned}
\end{equation}
  and the gradient of the potentail energy reads:
\begin{align}
  \nabla_{\tilde{\x}_0}\widetilde{\J}(\tilde{\x}_0)  &= -\left( \V^T \B_0 \V \right)^{-1} \left( \V^T\xb_0 - \tilde{\x}_0 \right) - \tilde{\boldsymbol{\lambdaup}}_0\,,
\end{align}
\end{subequations}
%
where $\tilde{\boldsymbol{\lambdaup}}_0$ is the solution of the ARRA reduced adjoint model \eqref{eqn:reduced_cost_function_gradient}.

The momentum variable is defined in the reduced space $\tilde{\p}_0 \in \mathbb{R}^{\nred}$, and the Hamiltonian reads:
\begin{equation}
\label{eqn:Reduced_Hamiltonian}
\widetilde{H}(\tilde{\p_0},\, \tilde{\x}_0) = \frac{1}{2} \, \tilde{\p}_0^T\, \widetilde{\M}^{-1}\, \tilde{\p}_0 + \widetilde{\J}(\tilde{\x}_0).
\end{equation}
Following \cite{Attia_HMCFilter_TR,Attia_HMCSmoother_TR}, the mass matrix  can be chosen as the diagonal matrix  $\widetilde{\M}={\rm diag}{(\V^T\B_0\V)^{-1}}\,$. Note that no further approximations are introduced to the numerical flow produced by the symplectic integrator because all calculations involving models states are calculated in the reduced space.

\subsection{Sampling in the full space using approximate gradients}

Here the model states are initially in the projected subspace defined by $\Pr=\V\V^T$ while the momentum is kept in the full space. The hope is that since the synthetic momentum is drawn at random from the full space for each proposed state, the symplectic integrator will help the sampler jump between slices of the full space rather that sampling a single subspace, leading to a better ensemble of states obtained from the original target posterior.

The target posterior distribution and the potential energy and its gradient are given by
\begin{subequations}
\label{eqn:Projected_HMC_Potential_Energy_and_Gradient}
\begin{equation}
\label{eqn:Full_PDF_with_Approximate_Gradient}
\begin{aligned}
\pi{(\xtin)} &= \exp{\big(- \wideparen{\J}(\xtin)\big)}\,, \\
  \wideparen{\J}\big(\xtin\big) &= \frac{1}{2} \lVert \wideparen{\x}_0 - \xb_0 \rVert^2 _{\B_0^{-1} } + \frac{1}{2}\sum_{k=0}^{\nobs}{ \lVert \y_k - \H_k(\wideparen{\x}_k)
                                    \rVert^2 _{ \mathbf{R}_k^{-1} } } \,, \\
  \wideparen{\x}_0 &= \xtin; \quad
  \wideparen{\x}_k = \V \widetilde{\mathcal{M}}_{k-1,\,k}(\V^T \wideparen{\x}_{k-1})\,,\quad  k=1,2,\ldots,\nobs.
\end{aligned}
\end{equation}
with gradient of the potential energy given by
\begin{align}
    \nabla_{\xtin}\, \wideparen{\J}\big(\xtin\big)  &= - {\mathbf{B}_0^{-1} ( \xb_0 - \xtin )} -\wideparen{\boldsymbol{\lambdaup}}_0,\,
\end{align}
\end{subequations}
where $\wideparen{\boldsymbol{\lambdaup}}_0$ is the solution of the following adjoint model
\begin{subequations}
\label{eqn:adjoint_double_proj}
\begin{eqnarray}
&& \qquad \wideparen{\boldsymbol{\lambdaup}}_{\nobs} = {\bf H}_{\nobs}^T \mathbf{R}_{\nobs}^{-1} \big({\bf y}_{\nobs} - \H_{\nobs}(\wideparen{\x}_{\nobs})\big), \quad \quad \quad \quad \quad \quad \quad \\
\nonumber
&& \qquad \wideparen{\boldsymbol{\lambdaup}}_{k-1} = \V\widetilde{\mathbf{M}}^T_{k-1,k}V^T\wideparen{\boldsymbol{\lambdaup}}_{k} + {\bf H}_{k-1}^T \mathbf{R}_{k-1}^{-1} \big({\mathbf y}_{k-1} - \H_{k-1}(\wideparen{\x}_{k-1})\big),\quad k=\nobs,..,1.
\end{eqnarray}
\end{subequations}
Here ${\bf H}_{k}$ represents the observation operator Jacobian linearized at $\wideparen{\x}_{k},~k=0,..,\nobs$, and $\widetilde{\mathbf{M}}_{k-1,k}$ is the Jacobian of the reduced order model evaluated at $V^T\wideparen{\x}_{k},~k=0,..,\nobs$.
The Hamiltonian in this case takes the form: 
\begin{equation}
\label{eqn:Projected_Hamiltonian}
{\wideparen H}(\p_0,\, \xtin) = \frac{1}{2} \, \p_0^T \mathbf{M}^{-1} \p_0 + {\wideparen \J}(\xtin)\,.
\end{equation}
An additional approximation is introduced to the numerical flow produced by the symplectic integrator by the approximation of the gradient of the potential energy. This may require more attention to be paid to the process of parameter tuning especially in the case of very high dimensional spaces.

Algorithm~\ref{alg:HMC_Sampling_Smoother_combined} summarizes the sampling process that yields an ensemble of states $\{\tilde{\x}_0(e) \,\in \mathbb{R}^{\nred}\}_{{e=1,2,\ldots,\nens}}$ in the reduced space, or ensemble of states $\{\wideparen{\x}_{0}(e),\,\in \mathbb{R}^{\nvar}\}_{{e=1,2,\ldots,\nens}}$ sampled from the high-fidelity state space with approximate gradient information, respectively
%

%
 \begin{algorithm}[htpb]
   \begin{algorithmic}[1]
 \State Initialize the mass matrix: $\widetilde{\M} \in \mathbf{R}^{\nred\times\nred}$ for sampling from~\eqref{eqn:Fully_Projected_PDF}, and $\M \in \mathbf{R}^{\nvar\times\nvar}$ for sampling from~\eqref{eqn:Full_PDF_with_Approximate_Gradient}.
 \State Initialize the chain. Preferably, the initial pair should be as close as possible to the target distribution.
 \State At each step $i$ of the Markov chain draw a random auxiliary momentum:
  $\tilde{\p}_0^{(i)} \sim \mathcal{N}(\mathbf{0}_{\nred},\widetilde{\M})$ for sampling from~\eqref{eqn:Fully_Projected_PDF}, and $\p_0^{(i)} \sim \mathcal{N}(\mathbf{0}_{\nred},\M)$ for sampling from~\eqref{eqn:Full_PDF_with_Approximate_Gradient}.
 \State Use a symplectic numerical integrator (e.g., position Verlet) to advance the current state by a pseudo-time increment $T$ to obtain a \textit{proposal}
  state :
\begin{equation}
\begin{aligned}
  \text{For sampling from~\eqref{eqn:Fully_Projected_PDF}: } ( \tilde{\p}_0^* , \tilde{\x}_0^* ) &= \widetilde{\phi}_T ( \tilde{\p}_0^{(i)} , {\tilde{\x}}_0^{(i)}  ). \\
  \text{For sampling from~\eqref{eqn:Full_PDF_with_Approximate_Gradient}: } ( {\p}_0^* , {\x}_0^* ) &= \wideparen{\phi}_T ( {\p}_0^{(i)} , {{\x}}_0^{(i)}  )\,,
\end{aligned}
\end{equation}
where $\wideparen{\Phi}_T$ indicates the flow approximation resulting from approximation of the gradient of the potential energy.
 \State For sampling from~\eqref{eqn:Fully_Projected_PDF}, use the Hamiltonian \eqref{eqn:Reduced_Hamiltonian} to evaluate the loss of energy`
  $\Delta \widetilde{H} = \widetilde{H}\big(\tilde{\p}_0^* , \tilde{\x}_0^*\big) - \widetilde{H}\big( \tilde{\p}_0^{(i)} , \tilde{\x}_0^{(i)}  \big) $.
  For sampling~\eqref{eqn:Full_PDF_with_Approximate_Gradient}, use the Hamiltonian~\eqref{eqn:Projected_Hamiltonian} to approximate the energy loss
  $\Delta \wideparen{H} = \wideparen{H}\big({\p}_0^* , {\x}_0^*\big) - \wideparen{H}\big( {\p}_0^{(i)} , {\x}_0^{(i)}  \big) $.
 \State Calculate the acceptance probability:
\begin{equation}
\label{eqn:acceptance_probability}
  \begin{aligned}
      \text{For sampling from~\eqref{eqn:Fully_Projected_PDF}: } a^{(i)} = 1 \wedge e^{-\Delta \widetilde{H}}\,, \\
      \text{For sampling from~\eqref{eqn:Full_PDF_with_Approximate_Gradient}: } a^{(i)} = 1 \wedge e^{-\Delta \wideparen{H}}\,.
  \end{aligned}
\end{equation}
 \State Discard both current and proposed momentum.
 \State \textbf{(Acceptance/Rejection)} Draw a uniform random variable $u^{(i)}\sim \mathcal{U}(0,1)$:
\begin{enumerate}
  \item[i-  ] If $a^{(i)} > u^{(i)}$ accept the proposal as the next sample; 
  \item[ii- ] If $a^{(i)} \leq u^{(i)}$ reject the proposal and continue with the  current state; 
\end{enumerate}
 \State Repeat steps $2$ to $7$ until $\nens$ distinct samples are drawn.
 \State Project the ensemble to the full space.
 \end{algorithmic}
 \caption{HMC Sampling ~\cite{Attia_HMCSmoother_TR}.}
  \label{alg:HMC_Sampling_Smoother_combined}
\end{algorithm}
%
Note that in Algorithm~\ref{alg:HMC_Sampling_Smoother_combined}, $\x_0^{(i)}$ refers to the model state at the initial time of the assimilation window (or models initial conditions) generated in step $i$ of the Markov chain.

\section{Properties of the Distributions Sampled with Reduced-Order Models}
\label{Sec:Reduced_Order_Sampling}
%
As explained above, our main goal in this work is to explore the possibility of lowering the computational expense posed by the original HMC smoother~\cite{Attia_HMCSmoother_TR} by following a reduced-order modeling approach. In the previous Section~\ref{Sec:Reduced_Order_Sampling},
we mensioned that the use of HMC sampling smoother with reduced order models requires following either of two alternatives, namely sampling the
posterior distribution fully projected in the lower dimensional subspace, or sampling the high fidelity distribution with gradients approximated
using information obtained from the reduced space.
In both cases, some amount of information will be lost due to either projecting the posterior PDF, or approximating the components appearing in the lielihood term.
More specifically, in the latter case, approximating the negative-log lielihood terms can lead to samples collected from a totaly different distribution than the
true posterior distribution.
In the rest of this section, we discuss the properties of the probability distributions resulting from projection or due to approximaton of the
negative-log likelihood terms making use of information coming only from a reduced-order subspace.

In the direct case where the posterior distribution is fully projected to the lower dimensional subspace, little can be said about the resulting distribution
unless if the true posterior is Gaussian. We explore this case in details in what follows.

\subsection{Projection of the posterior distribution for linear model and observation operators}
%
In this case the full distribution is projected into the lower-dimensional subspace by approximating both background and observation terms in Equation~\eqref{eqn:4DVar_Cost_Functional}. This projection leads to ensembles generated only in the reduced-space, and are then projected back to the high-fidelity space by left multiplication with $\V$. Projecting the ensembles back to the full space will not change their mass distribution  in the case of a linear model and observation operators, and will just embed the ensembles in the full space.

If both the model and the observation operator are linear operators, the posterior \eqref{eqn:Smoothing_Posterior_Kernel}  is a Gaussian distribution $\Pa(\xtin)= \N(\xa_0,\Aini)$, with a posterior (analysis) mean $\xa_0$, and an analysis error covariance matrix $\Aini$, i.e.
\begin{equation}
\label{eqn:Full_Posterior_PDF}
\Pa(\xtin) = \frac{(2 \pi)^{\frac{-\nvar}{2}}}{\sqrt{|\det(\Aini)|}} \exp{ \Bigl(- \frac{1}{2} \lVert \xtin-\xa_0 \rVert ^2 _{ \mathbf{\Aini}^{-1} } \Bigr) }\,.
\end{equation}
The mean and the covariance matrix of the Gaussian posterior~\eqref{eqn:Full_Posterior_PDF} are given by
\begin{equation}
\label{eqn:full_posterior_moments}
\begin{split}
\A_0^{-1} &= \B_0^{-1} + \sum_{k=0}^{\nobs} \M_{0,k}^T\, \mathbf{H}_k^T\, \R_k^{-1}\, \mathbf{H}_k\, \M_{0,k}\,, \\
\x_0^{\rm a} &= \A_0 \cdot \left( \B_0^{-1}\, \x_0^{\rm b} + \sum_{k=0}^{\nobs} \M_{0,k}^T\, \mathbf{H}_k^T\, \R_k^{-1}\, \y_k   \right).
\end{split}
\end{equation}

Projecting this PDF onto the subspace spanned by columns of the matrix $\V$ (e.g., POD basis) results in a projected PDF
$\tPa(\tilde{\x}_0) = \N(\V^T \xa_0, \V^T \Aini \V);\, \tilde{\x}_0 \in \mathbb{R}^{\nred}$, i.e.,
\begin{equation}
\label{eqn:POD_projected_Posterior}
\tPa(\tilde{\x}_0) = \frac{(2 \pi)^{\frac{-\nred}{2}}}{\sqrt{|\det{(\V^T \Aini \V)}| }}
\exp{ \Bigl(- \frac{1}{2} \lVert \tilde{\x}_0-\V^T\xa_0 \rVert ^2 _{\left( \mathbf{\V^T \Aini \V}\right)^{-1} } \Bigr) }\,.
\end{equation}
The linear transformation, of the analysis state, with the orthogonal projector $\Pr = \V \V^T$, results as well in the Gaussian distribution
$\dtPa(\wideparen{\x}_0 )= \N(\Pr \xa_0, \Pr \Aini \Pr) \equiv \N(\wideparen{\x}_0^\textnormal{a}, \wideparen{\A}_0)$, $\wideparen{\x}_0 \in \mathbb{R}^{\nvar}$.
The covariance matrix $\wideparen{\A}_0$ however is not full rank, and the Gaussian distribution is degenerate. The density function of this singular distribution can be  rigourously
formulated  by defining a restriction of Lebesgue measure to the affine subspace of $\mathbb{R}^{\nvar}$ whose dimension is limited to $\operatorname{rank}(\wideparen{\A}_0)$.
The Gaussian (singular) density then formula takes the form~\cite{khatri1968degenerateGaussian,rao2009linearinference}
\begin{equation}
\label{eqn:Singular_projected_Posterior}
\dtPa(\wideparen{\x}_0) = \frac{(2 \pi)^{\frac{-\nred}{2}}}{\sqrt{|\psdet{(\wideparen{\A}_0)}| }}
  \exp{ \Bigl(- \frac{1}{2} \lVert \wideparen{\x}_0 - \wideparen{\x}_0^\textnormal{a} \rVert ^2 _{  \wideparen{\A}_0 ^{\dagger} } \Bigr) }
                                         {\cdot \delta_{(\mathbf{I}-\Pr) \wideparen{\x}_0}}
\end{equation}
where $\psdet$ is the pseudo determinant, and $\dagger$ refers to the matrix pseudo inverse.
Of course, $\tilde{\x}_0 \in \mathbb{R}^{\nred},\, \wideparen{\x}_0 \in \mathbb{R}^{\nvar}$.
One can think of the PDF~\eqref{eqn:Singular_projected_Posterior} as a version of~\eqref{eqn:POD_projected_Posterior} embeded in the high-fidelity state space.
\begin{theorem}
If $\tPa(\tilde{\x}_0)$ and $\dtPa(\V {\bf \tilde x}_0)$ are the distributions defined in \eqref{eqn:POD_projected_Posterior} and \eqref{eqn:Singular_projected_Posterior}, respectively, then the following result holds true for a given reduced basis $\V$
\begin{equation}
\label{eqn:related_distributions}
\tPa(\tilde{\x}_0) = \dtPa( \V {\bf \tilde x}_0),\, \forall \tilde{\x}_0 \in \mathbb{R}^{\nred}.
\end{equation}
\end{theorem}
\begin{pf}
For this purpose it is sufficient to prove that
\begin{equation}
\label{eqn:Equivalence_Of_Projected_PDFs_RELATION}
  \lVert \wideparen{\x}_0- \wideparen{\x}_0^\textnormal{a} \rVert ^2 _{ \wideparen{\A}_0^{\dagger} }
   = \lVert \V^T\xtin-\V^T\xa_0 \rVert ^2 _{\left( \mathbf{\V^T \Aini \V}\right)^{-1} }
\end{equation}
Assume the relation given by Equation~\eqref{eqn:Equivalence_Of_Projected_PDFs_RELATION} is correct, we get the following:
\begin{subequations}
\label{eqn:Equivalence_Of_Projected_PDFs_PROOF}
 \begin{equation}
 \begin{aligned}
\lVert \wideparen{\x}_0-\wideparen{\x}_0^\textnormal{a} \rVert ^2 _{ {\wideparen{\A}_0}^{\dagger} }
 &= \lVert \V^T \xtin-\V^T\xa_0 \rVert ^2 _{\left( {\V^T \Aini \V}\right)^{-1} }\,,  \\
(\wideparen{\x}_0-\wideparen{\x}_0^\textnormal{a})^T { \wideparen{\A}_0^{\dagger} } (\wideparen{\x}_0-\wideparen{\x}_0^\textnormal{a})
 &= (\V^T \xtin-\V^T\xa_0)^T {\left( {\V^T \Aini \V}\right)^{-1} } (\V^T \xtin-\V^T\xa_0)  \\
(\Pr\xtin-\Pr\xa_0)^T { \wideparen{\A}_0^{\dagger} } (\Pr\xtin-\Pr\xa_0)
 &= (\V^T \xtin-\V^T\xa_0)^T {\left( {\V^T \Aini \V}\right)^{-1} } (\V^T \xtin-\V^T\xa_0)
 \end{aligned}
 \end{equation}
 Or equivalently:
 \begin{equation}
 \begin{aligned}
0 &= \left(\Pr\xtin-\Pr\xa_0 \right)^T { \wideparen{\A}_0^{\dagger} } \left(\Pr\xtin-\Pr\xa_0 \right)
 - (\V^T \xtin-\V^T\xa_0)^T {\left( {\V^T \Aini \V}\right)^{-1} } (\V^T \xtin-\V^T\xa_0)  \,  \\
  &= \left(\V ( \V^T\xtin-\V^T\xa_0) \right)^T { \wideparen{\A}_0^{\dagger} } \left(\V ( \V^T\xtin-\V^T\xa_0) \right)
 - \left(\V^T \xtin-\V^T\xa_0 \right)^T {\left( {\V^T \Aini \V}\right)^{-1} } \left(\V^T \xtin-\V^T\xa_0\right)\,  \\
  &= \left( \V^T\xtin-\V^T\xa_0\right)^T \V^T{ \wideparen{\A}_0^{\dagger} } \V \left( \V^T\xtin-\V^T\xa_0 \right)
 - \left(\V^T \xtin-\V^T\xa_0 \right)^T {\left( {\V^T \Aini \V}\right)^{-1} } \left(\V^T \xtin-\V^T\xa_0 \right)\,  \\
  &= ( \V^T\xtin-\V^T\xa_0)^T \left(
  {\V^T{\wideparen{\A}_0^{\dagger} } \V }-
  {\left( {\V^T \Aini \V}\right)^{-1} }
 \right) ( \V^T\xtin-\V^T\xa_0) \,.
 \end{aligned}
 \end{equation}
This holds true if the matrix ${\V^T{ \wideparen{\A}_0^{\dagger} } \V }- {\left( {\V^T \Aini \V}\right)^{-1} }$ is equal to a zero matrix.
The matrix $\V$ has orthonormal columns and consequently $\left( {\Pr \Aini \Pr}\right)^{\dagger} = \left( {\V \V^T \Aini \Pr}\right)^{\dagger} = \left( {\V^T \Aini \Pr}\right)^{\dagger} \V^{\dagger}$.
Since the pseudo inverse and the transpose operations are commutative, we get the following:
 \begin{equation}
 \begin{aligned}
\left( {\V^T \Aini \Pr}\right)^{\dagger} &= \left( \left( { \Pr \Aini^T \V}\right)^T \right) ^{\dagger}\,,  \\
 &= \left(\V^T \right)^{\dagger}  \left( { \V^T \Aini \V}\right)^{\dagger} \,,
 \end{aligned}
 \end{equation}
 and consequently:
 \begin{equation}
 \begin{aligned}
{\V^T{ \wideparen{\A}_0^{\dagger} } \V } = {\V^T{\left( {\Pr \Aini \Pr}\right)^{\dagger} } \V }
 &= \V^T \left(\V^T \right)^{\dagger}  \left( { \V^T \Aini \V}\right)^{\dagger} \V^{\dagger} \V \,,  \\
 &= \V^T \left(\V^T \right)^{\dagger}  \left( { \V^T \Aini \V}\right)^{-1} \V^{\dagger} \V \,,  \\
 &= (\V^T \V)  \left( { \V^T \Aini \V}\right)^{-1} (\V^T \V) \,,  \\
 &= \left( { \V^T \Aini \V}\right)^{-1} \,,
 \end{aligned}
 \end{equation}
\end{subequations}
where $\V^{\dagger} = \V^T,\, \text{and} \left(\V^T\right)^{\dagger} = \V$ since $\V$ has orthonormal columns.
This means that the relation ~\eqref{eqn:Equivalence_Of_Projected_PDFs_RELATION} holds,
and the equivalence between~\eqref{eqn:POD_projected_Posterior} and~\eqref{eqn:Singular_projected_Posterior} follows immediately.
\end{pf}
This result suggests that sampling from the distribution~\eqref{eqn:Singular_projected_Posterior} can be carried out efficiently by sampling the distribution~\eqref{eqn:POD_projected_Posterior}, then projecting the ensembles back to the full space using $\V$.

By determining the Kullback Leibler (KL)~\cite{cover2012elementsKLD} divergence measure between the high fidelity distribution $\Pa(\xtin)$ and the probability distribution $\dtPa(\wideparen{\x}_0)$, one can estimate the error between the projected samples obtained using distribution \eqref{eqn:POD_projected_Posterior} and those sampled from the high fidelity distribution $\Pa(\xtin)$.
\begin{theorem}
The KL divergence measure between the Gaussian distribution $\dtPa(\xtin)$ given by \eqref{eqn:Singular_projected_Posterior}, and the probability distribution $\Pa(\xtin)$ defined in~\eqref{eqn:Full_Posterior_PDF}, is given as
\begin{equation}
\label{eqn:thoerem_KL_measure_Full_vs_projected}
\begin{split}
   D_{\rm KL} \left( \dtPa(\xtin) || \Pa(\xtin) \right)
   &= \frac{1}{2} \left( (\nvar-\nred) \ln{ (2 \pi)}
                   + \ln\left( \frac{|\det(\Aini)|}{|\psdet{(\wideparen{\A}_0)}|} \right)
                   + \lVert \wideparen{\x}_0^\textnormal{a} - \xa_0 \rVert _ {\A_0^{-1}}
                   + \textnormal{trace}{
                                        \left(
                                              \left( \A_0^{-1} -\wideparen{\A}_0^{\dagger}\right)\, \wideparen{\A}_0
                                        \right)
                                        }
                  \right)
\,,
\end{split}
\end{equation}
where $\V \in \mathbf{R}^{\nvar \times \nred}$ and $\nred < \nvar$.
\end{theorem}
\begin{pf}
The KL measure is obtained as
\begin{subequations}
\label{eqn:KL_Divergence_1}
\begin{align}
D_{\rm KL} \left( \dtPa(\xtin) || \Pa(\xtin) \right) &=
\mathtt{E}_{\dtPa} \left[ \ln \left( \frac{\dtPa(\xtin)}{\Pa(\xtin)} \right)  \right]\,, \\
\ln \left( \frac{\dtPa(\xtin)}{\Pa(\xtin)} \right)
  &=   \ln{\left(
\frac{(2 \pi)^{\frac{-\nred}{2}}}{\sqrt{|\psdet{(\wideparen{\A}_0)}| }}
  \right)}
  + \ln{\left(
\frac{\sqrt{|\det(\Aini)|}}{(2 \pi)^{\frac{-\nvar}{2}} }
  \right)}
  +
\ln{\left(
\frac{\exp{ \Bigl(- \frac{1}{2} \lVert \xtin-\wideparen{\x}_0^\textnormal{a} \rVert ^2 _{ {\wideparen{\A}_0}^{\dagger} } \Bigr) }}
{\exp{ \Bigl(- \frac{1}{2} \lVert \xtin-\xa_0 \rVert ^2 _{ \mathbf{\Aini}^{-1} } \Bigr) }}
\right)} \\
  %
%
&= \frac{\left(\nvar-\nred\right) \ln{ (2 \pi)} }{2}  + \frac{1}{2}\ln\left( \frac{|\det(\Aini)|}{|\psdet{(\wideparen{\A}_0)}|} \right)
  + \frac{1}{2} {   \left( \lVert \xtin-\xa_0 \rVert ^2 _{ \mathbf{\Aini}^{-1} }
  - \lVert \xtin-\wideparen{\x}_0^\textnormal{a} \rVert ^2 _{ {\wideparen{\A}_0}^{\dagger} } \right)
}  \\
\mathtt{E}_{\dtPa}\left[ \ln \left( \frac{\dtPa(\xtin)}{\Pa(\xtin)} \right)  \right]
&= \frac{\left(\nvar-\nred\right) \ln{ (2 \pi)} }{2}  + \frac{1}{2} \ln\left( \frac{|\det(\Aini)|}{|\psdet{(\wideparen{\A}_0)}|} \right)
  +\frac{1}{2} \mathtt{E}_{\dtPa}\left[
  \left( \lVert \xtin-\xa_0 \rVert ^2 _{ \mathbf{\Aini}^{-1} }
- \lVert \xtin-\wideparen{\x}_0^\textnormal{a} \rVert ^2 _{ {\wideparen{\A}_0}^{\dagger} } \right)
\right]\,, \label{eqn:KL_Divergence_Expectation}
\end{align}
\end{subequations}
where $\ln{\left( \frac{|\det(\Aini)|}{|\psdet{(\wideparen{\A}_0)}|} \right)}$ is the sum of logarithms of eigenvalues of $\Aini$ lost due to projection. This value can be also replaced with
$ln{\left(\frac{|\det(\Aini)|}{|\det{(\V^T \Aini \V)}|}\right)}$ due to the nature of the matrix $\V$. The expectation of the quadratic terms in
Equation~\eqref{eqn:KL_Divergence_Expectation}
can be obtained as follows:
\begin{subequations}
\label{eqn:Quadratic_Terms_KL_Divergence}
\begin{align}
\mathtt{E}_{\dtPa}\left[
\left( \lVert \xtin-\xa_0 \rVert ^2 _{ {\Aini}^{-1} }
  - \lVert \xtin-\wideparen{\x}_0^\textnormal{a} \rVert ^2 _{ {\wideparen{\A}_0}^{\dagger} } \right)
  \right]
&= \mathtt{E}_{\dtPa}\left[
\left( \lVert \xtin-\xa_0 \rVert ^2 _{ {\Aini}^{-1} } \right]
  -
  \mathtt{E}_{\dtPa}\left[
  \lVert \xtin-\wideparen{\x}_0^\textnormal{a} \rVert ^2 _{ {\wideparen{\A}_0}^{\dagger} } \right) \right]\\
&= \mathtt{E}_{\dtPa}\left[
( \xtin-\xa_0 )^T {\Aini}^{-1} ( \xtin-\xa_0 )
\right] -
  \mathtt{E}_{\dtPa}\left[
(\xtin-\wideparen{\x}_0^\textnormal{a} )^T  {\wideparen{\A}_0}^{\dagger}  (\xtin-\wideparen{\x}_0^\textnormal{a} )
\right] \\
&= (\wideparen{\x}_0^\textnormal{a} - \xa_0)^T \A_0^{-1} (\wideparen{\x}_0^\textnormal{a} - \xa_0) +
  \tr{ \left( \A_0^{-1} \, {\wideparen{\A}_0} \right) }
  - \tr{ \left( {\wideparen{\A}_0}^{\dagger} \, {\wideparen{\A}_0} \right) }
\end{align}
\end{subequations}
from Equations~\eqref{eqn:Quadratic_Terms_KL_Divergence}, and~\eqref{eqn:KL_Divergence_1}, we obtain:

\begin{subequations}
\begin{align}
\label{eqn:KL_Divergence_2}
D_{\rm KL} \left( \dtPa(\xtin) || \Pa(\x_0) \right) &=
  \frac{\left(\nvar-\nred\right) \ln{ (2 \pi)} }{2}  + \frac{1}{2} \ln\left( \frac{|\det(\Aini)|}{|\psdet{(\wideparen{\A}_0)}|} \right)
+ \frac{1}{2} \lVert \wideparen{\x}_0^\textnormal{a} - \xa_0 \rVert _ {\A_0^{-1}} \\
  &\quad +
  \frac{1}{2} \tr{ \left( \A_0^{-1} \, {\wideparen{\A}_0} \right) }
  - \frac{1}{2} \tr{ \left( {\wideparen{\A}_0}^{\dagger}  \, {\wideparen{\A}_0} \right) }   \nonumber  \\
  &= \frac{\left(\nvar-\nred\right) \ln{ (2 \pi)} }{2}  + \frac{1}{2} \ln\left( \frac{|\det(\Aini)|}{|\psdet{(\wideparen{\A}_0)}|} \right)
+ \frac{1}{2} \lVert \wideparen{\x}_0^\textnormal{a} - \xa_0 \rVert _ {\A_0^{-1}} \\
  &\quad +
  \frac{1}{2} \tr{
\left(
\A_0^{-1} \, {\wideparen{\A}_0} - {\wideparen{\A}_0}^{\dagger} \, {\wideparen{\A}_0}
\right)
}   \ \nonumber  \\
  &= \frac{1}{2} \left( \left(\nvar-\nred\right) \ln{ (2 \pi)}
+  \ln\left( \frac{|\det(\Aini)|}{|\psdet{(\wideparen{\A}_0)}|} \right)  +   \lVert \wideparen{\x}_0^\textnormal{a} - \xa_0 \rVert _ {\A_0^{-1}}
+ \textnormal{trace}{
\left(
\left( \A_0^{-1} -\wideparen{\A}_0^{\dagger}\right)\, \wideparen{\A}_0
\right)
}
\right)
\,,
\end{align}
\end{subequations}
which completes the proof.
\end{pf}
This measure can be used to quantify the quality of POD basis given an estimation of the analysis error covariance matrix, e.g., based on an ensemble of states, sampled from the high fidelity distribution, or approximated based on statistics of the 4D-Var cost functional.
Notice that the KL measure given in~\eqref{eqn:thoerem_KL_measure_Full_vs_projected} is finite since $\dtPa(\xtin)$ is absolutely continuous with respect to $\Pa(\xtin)$ (and it is zero only if $\dtPa(\xtin)=\Pa(\xtin)$).
For this reason, we set the projected PDF as the reference density in the KL measure.

\subsection{Approximating the likelihood function using reduced order models}
In the latter approach, the background term is kept in the high fildelity space, while only the terms involving model propagations are approximated using reduced-order models. This means that the target distribution is the PDF give by~\eqref{eqn:Full_PDF_with_Approximate_Gradient}.
The use of this approximation in the HMC algorithm results in samples collected from the distribution~\eqref{eqn:Full_PDF_with_Approximate_Gradient}.
This approximation maintains the background term in the full space, while the model states involved in the observation term are approximated in the lower-dimensional subspace. This means that the posterior distribution is non-degenerate in the full space due to the background term. However, it is not immediately obvious which distributions samples will be collected from. In Theorem~\ref{tho:Approximate_Smoothing_Posterior_Kernel_1_update} we show the link between
posterior distribution given by~\eqref{eqn:Full_PDF_with_Approximate_Gradient} and the distribution defined in \eqref{eqn:Smoothing_Posterior_Kernel}.
\begin{theorem}
\label{tho:Approximate_Smoothing_Posterior_Kernel_1_update}
  The posterior distribution ${\boldsymbol{\pi}}$ defined in \eqref{eqn:Smoothing_Posterior_Kernel} associated with the high-fidelity model \eqref{eqn:Forward_Model} is proportional to the analysis posterior distribution $\widetilde{\boldsymbol{\pi}}$ introduced in \eqref{eqn:Full_PDF_with_Approximate_Gradient} associated with the reduced order model, by the ratio between joint likelihood functions given projected and high-fidelity states, i.e.
  \begin{equation}
 \widetilde{\boldsymbol{\pi}}(\xtin)
   =  \boldsymbol{\pi}(\xtin) \cdot
\prod\limits_{k=0}^{\nobs} \frac{ \mathcal{P}(\y_k |  \x_k = \V\, \widetilde{\mathcal{M}}_{0,k}(\V^T\x_0) \big) }
{ \mathcal{P}\big(\y_k | \x_k = \mathcal{M}_{0,k}(\x_0) \big) }\,.
  \end{equation}

\end{theorem}

\begin{pf}
  The exact and the approximate posterior distributions $\boldsymbol{\pi}(\xtin),\, \widetilde{\boldsymbol{\pi}}(\xtin) $ are generally described as follows
  \begin{subequations}
  \label{eqn:true_and_approximate_posterior_probabilities}
  \begin{equation}
  \begin{split}
\boldsymbol{\pi}(\xtin) = \Pa(\xtin) &= \Pb(\xtin) \cdot  \mathcal{P}(\y_0 | \x_0) \cdot \prod\limits_{k=1}^{\nobs} \mathcal{P}(\y_k | \x_k) \cdot \mathcal{P}(\x_k | \x_{k-1}) \\
&= \Pb(\xtin) \cdot  \mathcal{P}(\y_0 | \x_0) \cdot \prod\limits_{k=1}^{\nobs} \mathcal{P}\big(\y_k | \x_k = \mathcal{M}_{0,k}(\x_0) \big)\,,  \\
  \end{split}
  \end{equation}
  \begin{equation}
  \begin{split}
\widetilde{\boldsymbol{\pi}}(\xtin) &= \Pb(\xtin) \cdot  \mathcal{P}(\y_0 | \tilde{\x}_0) \cdot \prod\limits_{k=1}^{\nobs}
  \mathcal{P}(\y_k | \tilde{\x}_k) \cdot \mathcal{P}(\tilde{\x}_k | \tilde{\x}_{k-1}) \\
&= \Pb(\xtin) \cdot  \mathcal{P}(\y_0 | \tilde{\x}_0) \cdot \prod\limits_{k=1}^{\nobs} \mathcal{P}(\y_k |  \x_k = \V\, \widetilde{\mathcal{M}}_{0,k}(\V^T\x_0) \big).
  \end{split}
  \end{equation}
  This leads to the following:
  \begin{equation}
  \begin{split}
\frac{\widetilde{\boldsymbol{\pi}}(\xtin)}{\boldsymbol{\pi}(\xtin)}
   &= \frac{ \prod\limits_{k=0}^{\nobs} \mathcal{P}(\y_k |  \x_k = \V\, \widetilde{\mathcal{M}}_{0,k}(\V^T\x_0) \big) }
   { \prod\limits_{k=0}^{\nobs} \mathcal{P}\big(\y_k | \x_k = \mathcal{M}_{0,k}(\x_0) \big) } \,, \\
\widetilde{\boldsymbol{\pi}}(\xtin)
   &=  \boldsymbol{\pi}(\xtin) \cdot
\prod\limits_{k=0}^{\nobs} \frac{ \mathcal{P}(\y_k |  \x_k = \V\, \widetilde{\mathcal{M}}_{0,k}(\V^T\x_0) \big) }
{ \mathcal{P}\big(\y_k | \x_k = \mathcal{M}_{0,k}(\x_0) \big) }\,. \\
  \end{split}
  \end{equation}
  \end{subequations}
\end{pf}

This result suggests that the larger the distances $\| \x_k - \V {\tilde x}_{k}\|_2,~k=1,2,..,\nobs$ are, the more different the distributions ${\boldsymbol{\pi}}$ and $\widetilde{\boldsymbol{\pi}}$ will be. By selecting appropriate reduced manifolds $\V$ and decreasing the error associated with the reduced order models, the ratio can be brought closer to $1$.

\begin{corollary}
The KL divergence measure between the original posterior~\eqref{eqn:Smoothing_Posterior_Kernel} and the approximated distribution \eqref{eqn:Full_PDF_with_Approximate_Gradient}
is:
\begin{equation}
\begin{aligned}
D_{\rm KL}(\widetilde{\boldsymbol{\pi}} ||\boldsymbol{\pi})
&= \mathtt{E}_{\widetilde{\boldsymbol{\pi}}}\left[ \ln ( \widetilde{\boldsymbol{\pi}} ) - \ln (\boldsymbol{\pi}) \right] = \mathtt{E}_{\widetilde{\boldsymbol{\pi}}}\left[ \J(\xtin) -\widetilde{\J}(\xtin) \right]\\
&= \mathtt{E}_{\widetilde{\boldsymbol{\pi}}}\left[ \J^\textnormal{obs}(\xtin) -\widetilde{\J}^\textnormal{obs}(\xtin) \right]\,,
\end{aligned}
\end{equation}
where $\J^\textnormal{obs}(\xtin),\, \text{and } \widetilde{\J}^\textnormal{obs}(\xtin)$ are the observation terms in the full and the approximate 4D-Var cost function.
\end{corollary}

\begin{corollary} In the filtering case, where only one observation is assimilated, if the initial condition is projected on the columns of $\V$ to approximate the likelihood term, the posterior distribution is given by:
  \begin{equation}
\widetilde{\boldsymbol{\pi}}(\xtin) = \widetilde{\Pa}(\xtin)
\propto \frac{\boldsymbol{\pi}(\xtin^\parallel)\cdot \Pb(\xtin)}{\Pb(\xtin^\parallel)}\,,
  \end{equation}
  where $\xtin = \xtin^\parallel + \xtin^\perp$, with $\x_0^\parallel \in range(\V)$ and $\x_0^\perp \in null(\V^T)$
\end{corollary}

\begin{pf}
The two states $\xtin$ and $\xtin^\parallel$ differ only along a direction $\xtin^\perp$ orthogonal to the reduced space, that is $\V^T \xtin^\perp = 0$, and consequently
\[
\V^T\, \xtin = \V^T(\xtin^\parallel + \xtin^\perp) = \V^T\, \xtin^\parallel.
\]
  In the filtering case the cost function reads:
  \begin{subequations}
  \label{eqn:projected_filtering_posterior_2}
  \begin{equation}
  \begin{split}
 \widetilde{\J}(\xtin) &= \frac{1}{2} \lVert \xtin^\parallel + \xtin^\perp - \xb_0 \rVert ^2 _{ \Bini^{-1} }
   + \frac{1}{2} { \left\lVert \y_k - \H_0\left(\V{ \V^T \xtin^\parallel }\right) \right\rVert _{ \mathbf{R}_0^{-1} } } \\
   &= \frac{1}{2} \lVert \xtin^\parallel + \xtin^\perp - \xb_0 \rVert ^2 _{ \Bini^{-1} }
   + \frac{1}{2}{ \left\lVert \y_k - \H_0 ( \xtin^\parallel)  \right\rVert _{ \mathbf{R}_0^{-1} } } \,, \\
{\J}(\xtin^\parallel ) &= \frac{1}{2} \lVert \xtin^\parallel - \xb_0 \rVert ^2 _{ \Bini^{-1} }
 + \frac{1}{2}{ \left\lVert \y_k - \H_0\left({ \xtin^\parallel }\right) \right\rVert _{ \mathbf{R}_0^{-1} } }, \\
 - \widetilde{\J}(\xtin) &= - \J(\xtin^\parallel )
   - \frac{1}{2} \lVert \xtin^\parallel + \xtin^\perp - \xb_0 \rVert ^2 _{ \Bini^{-1} }
   + \frac{1}{2} \lVert \xtin^\parallel - \xb_0 \rVert ^2 _{ \Bini^{-1} }.
  \end{split}
  \end{equation}
  \end{subequations}
   Exponentiating of both sides leads to the following:
  \begin{equation}
  \begin{aligned}
 \exp{\left( - \widetilde{\J}(\xtin) \right)} &= \exp{ \left( -\frac{1}{2} \lVert \xtin^\parallel + \xtin^\perp - \xb_0 \rVert ^2 _{ \Bini^{-1} } \right)}
   \exp{\left( - \J(\xtin^\parallel ) \right)}
   \exp{\left( \frac{1}{2} \lVert \xtin^\parallel - \xb_0 \rVert ^2 _{ \Bini^{-1} } \right)} \,,  \\
 {\widetilde{\boldsymbol{\pi}}}(\xtin)  &\propto  \frac{\boldsymbol{\pi}(\xtin^\parallel)\cdot \Pb(\xtin)}{\Pb(\xtin^\parallel)}\,.
  \end{aligned}
  \end{equation}
  This completes the proof. 
\end{pf}

\begin{corollary}
  In the filtering case, if $\x_0^\perp=0$ the two distributions $\widetilde{\boldsymbol{\pi}}(\xtin)$, and $\boldsymbol{\pi}(\xtin)$ coincide, and
  if $\x_0^\parallel=0$ then the reduced distribution $\widetilde{\boldsymbol{\pi}}(\xtin)$ coincides with the background distribution $\Pb(\xtin)$.
\end{corollary}

In the general case we have:
  \begin{subequations}
  \label{eqn:projected_filtering_posterior_1}
  \begin{equation}
  \begin{aligned}
 \widetilde{\J}(\xtin) = \widetilde{\J}(\xtin^\parallel + \xtin^\perp)
   &= \frac{1}{2} \lVert \xtin^\parallel + \xtin^\perp - \xb_0 \rVert ^2 _{ \Bini^{-1} }
   + \frac{1}{2}\sum_{k=0}^{\nobs}{ \left\lVert \y_k - \H_k\left(\V{ {\widetilde{\mathcal{M}}_{0,k}}\left(\V^T (\xtin^\parallel + \xtin^\perp)\right) }\right) \right\rVert _{ \mathbf{R}_k^{-1} } } \\
   &= \frac{1}{2} \lVert \xtin^\parallel + \xtin^\perp - \xb_0 \rVert ^2 _{ \Bini^{-1} }
   + \frac{1}{2}\sum_{k=0}^{\nobs}{ \left\lVert \y_k - \H_k\left(\V{ {\widetilde{\mathcal{M}}_{0,k}}\left(\V^T \xtin^\parallel\right) }\right) \right\rVert _{ \mathbf{R}_k^{-1} } } \,,
  \end{aligned}
  \end{equation}
  \begin{equation}
 {\J}(\xtin^\parallel ) = \frac{1}{2} \lVert \xtin^\parallel - \xb_0 \rVert ^2 _{ \Bini^{-1} }
 + \frac{1}{2}\sum_{k=0}^{\nobs}{ \left\lVert \y_k - \H_k\left({ {{\mathcal{M}}_{0,k}}\left(\xtin^\parallel\right) }\right) \right\rVert _{ \mathbf{R}_k^{-1} } } \,,
  \end{equation}
  \end{subequations}

\begin{corollary}
The posterior $\widetilde{\boldsymbol{\pi}}$  \eqref{eqn:Full_PDF_with_Approximate_Gradient} is Gaussian with analysis covariance and mean:
\begin{equation}
\label{eqn:approximate_posterior_moments}
\begin{split}
{\wideparen{\A}}_0^{-1} &= \B_0^{-1} + \sum_{k=0}^{\nobs} \V\, \widetilde{\M}_{0,k}^T\, \mathbf{H}_k^T\, \R_k^{-1}\, \mathbf{H}_k\, \widetilde{\M}_{0,k}\,\V^T \\
{\wideparen{\x}}_0^{\rm a} &= {\wideparen{\A}}_0 \cdot \left( \B_0^{-1}\, \x_0^{\rm b} + \sum_{k=0}^{\nobs} \V \widetilde{\M}_{0,k}^T\, \mathbf{H}_k^T\, \R_k^{-1}\, \y_k   \right).
\end{split}
\end{equation}
\end{corollary}
From Equations~\eqref{eqn:full_posterior_moments} and~\eqref{eqn:approximate_posterior_moments} we conclude that the analysis mean and covariance associated with the distribution $\widetilde{\boldsymbol{\pi}}$ \eqref{eqn:Full_PDF_with_Approximate_Gradient} are not obtained simply by projecting the mean and covariance of the high fidelity distribution ${\boldsymbol{\pi}}$ \eqref{eqn:Smoothing_Posterior_Kernel}, i.e., ${\wideparen{\A}}_0 \neq \V\V^T \A_0 \V\V^T$ and ${\wideparen{\x}}_0^{\rm a} \neq \V\V^T \x_0^{\rm a}$.

\begin{corollary}
For a constant model operator $\M_{k-1,k}=\M$ the mean and the covariance of the high-fidelity posterior \eqref{eqn:Smoothing_Posterior_Kernel} are
\begin{equation}
\label{eqn:full_posterior_moments_linear_model}
\begin{split}
\A_0^{-1} &= \B_0^{-1} + \sum_{k=0}^{\nobs} (\M^k)^T\, \mathbf{H}_k^T\, \R_k^{-1}\, \mathbf{H}_k\, (\M^k)\,, \\
\x_0^{\rm a} &= \A_0 \cdot \left( \B_0^{-1}\, \x_0^{\rm b} + \sum_{k=0}^{\nobs} (\M^k)^T\, \mathbf{H}_k^T\, \R_k^{-1}\, \y_k   \right)\,,
\end{split}
\end{equation}
while in the case of posterior \eqref{eqn:Full_PDF_with_Approximate_Gradient}, the associated analysis covariance and mean are
\begin{equation}
\label{eqn:approximate_posterior_moments_linear_model}
\begin{split}
{\wideparen{\A}}_0^{-1} &= \B_0^{-1} + \sum_{k=0}^{\nobs} \left( (\Pr \M \Pr)^k \right)^T\, \mathbf{H}_k^T\, \R_k^{-1}\, \mathbf{H}_k\, (\Pr \M \Pr)^k \\
{\wideparen{\x}}_0^{\rm a} &= {\wideparen{\A}}_0 \cdot \left( \B_0^{-1}\, \x_0^{\rm b} + \sum_{k=0}^{\nobs} \left( (\Pr \M \Pr)^k \right)^T\, \mathbf{H}_k^T\, \R_k^{-1}\, \y_k   \right).
\end{split}
\end{equation}
\end{corollary}
A closed form for the distribution~\eqref{eqn:Full_PDF_with_Approximate_Gradient} can be obtained if
a) the observation errors are defined given the state vectors in the lower-dimensional subspace embedded in the full space (the projected space),
b) the observation errors at time $t_k$ follow Gaussian distribution with zero mean and covariance matrix $\mathbf{\tilde{R}}_k$, that is
\begin{equation}
\label{eqn:Projected_Observation_Errors}
  \mathbf{\tilde{e}}_k^{\rm obs} = \yk - \H(\V\tilde{\x}_k) \sim \N(0,\, \mathbf{\tilde{R}}_k)\,,
\end{equation}
and c) forcing the regular assumptions of time independence of observation errors, and independence from model background state (in the smaller space), once can obtain the posterior defined by~\eqref{eqn:Full_PDF_with_Approximate_Gradient}.

\section{Numerical Results}
\label{Sec:Numerical_Results}
In this section we test numerically the reduced order sampling algorithms using the shallow-water equations (SWE) model in  Cartesian coordinates.
%
\subsection{The SWE model}
%
Many phenomena in fluid dynamics are characterized by horizontal length scale much greater than the vertical length, consequently when equipped with Coriolis forces,
the shallow water equations model (SWE) becomes a valuable tool in atmospheric modeling, as a simplification of the primitive equations of atmospheric flow.
Their solutions represent many of the types of motion found in the real atmosphere, including slow-moving Rossby waves and fast-moving gravity waves~\cite{heikes1995numerical}.
The alternating direction fully implicit finite difference scheme~\cite{Gus1971} was considered in this paper and it is stable for large CFL condition numbers.

The SWE model using the $\beta$-plane approximation on a rectangular domain is introduced (see~\cite{Gus1971})
\begin{equation}\label{eqn:swe-pde}
\frac{\partial w}{\partial t}=A(w)\frac{\partial w}{\partial x}+B(w)\frac{\partial w}{\partial y}+C(y)w,
\quad (x,y) \in [0,L] \times [0,D], \quad t\in(0,t_{\rm f}],
\end{equation}
where $w=(u,v,\phi)^T$ is a vector function, $u,v$ are the velocity components in the $x$ and $y$ directions, respectively, $h$ is the depth of the fluid, $g$ is the acceleration
due to gravity, and $\phi = 2\sqrt{gh}$.

The matrices $A$, $B$ and $C$ have the form
\begin{equation}
A=-\left[\begin{array}{ccc}
   u&0&\phi/2\\
   0&u&0\\
   \phi/2 &0&u
 \end{array}\right],
 \quad
B=-\left[\begin{array}{ccc}
   v&0&0\\
   0&v&\phi/2\\
   0&\phi/2 &v
 \end{array}\right],
 \quad
C=\left[\begin{array}{rrr}
   0&f&0\\
   -f&0&0\\
   0&0&0
  \end{array}\right]\,,
\end{equation}
where $f$ is the Coriolis term
\begin{equation}
f=\hat f + \beta(y-D/2),\quad \beta=\frac{\partial f}{\partial y}, \quad \forall\, y \in [0,D],
\end{equation}
with $\hat f$ and $\beta$ constants.
We assume periodic solutions in the $x$ direction for all three state variables
while in the $y$ direction $v(x,0,t)=v(x,D,t)=0,\quad x\in[0,L],\quad t\in(0,t_{\rm f}]\,,\,$
and Neumann boundary condition are considered for $u$ and $\phi$.

The numerical scheme is implemented in Fortran and uses a sparse matrix environment. For operations with sparse matrices we utilize SPARSEKIT library \cite{Saad1994},
and the sparse linear systems resulted from quasi-Newton iterations is solved using MGMRES library \cite{Barrett94,Kelley95,Saad2003}.
More details on the implementation can be found in \cite{cstefuanescu2015pod}.

\subsection{Smoothing experimental settings}

To test the HMC smoothers with SWE model in the context for data assimilation, we construct an assimilation window of length $91$ units, with $10$ observations distributed
over the window. Here the observations are linearly related to model state with $\H = \mathbf{I}$, where $\mathbf{I}$ is the identity matrix.
4D-Var is carried out in both high-fidelity space (Full 4D-Var) and reduced-order space (Reduced 4D-Var) against the HMC sampling smoother in the following settings
\begin{enumerate}
\item[i)  ] Sampling the high-fidelity space using the original HMC smoother~\cite{Attia_HMCSmoother_TR} (``Full HMC''),
\item[ii) ] Sampling the reduced space, i.e. sampling~\eqref{eqn:Fully_Projected_PDF} (``Reduced HMC''),
\item[iii)] Sampling the high-fidelity space with approximate gradients, i.e. sampling~\eqref{eqn:Full_PDF_with_Approximate_Gradient} (``Approximate Full HMC''),
\end{enumerate}

In the three cases, the symplectic integrator used is the position Verlet~\eqref{eqn:Position_Verlet} with step size parameters tuned empirically through a preprocessing step. Higher order integrators~\cite{Attia_HMCFilter_TR,sanz2014markov} and automatic tuning of parameters should be considered when theses algorithms are applied to more complicated,
e.g., when $\H$ is nonlinear or when the Gaussian prior assumption is relaxed. The reduced basis $\V$ is constructed using  initial trajectories of the high-fidelity forward and adjoint models as well as the associated gradient of the full cost function \cite{cstefuanescu2015pod}. Later on this basis is updated using the current proposal and the corresponding trajectories.
%
\subsection{Numerical results}

Due to the simple settings described above the posterior distribution is not expected to deviate notably from a Gaussian. This will enable us to easily test the quality of the ensemble by testing the first two moments generated from the ensemble.

\begin{figure}
\centering
\caption{Data assimilation results using 4D-Var schemes, and HMC smoother, in both high-fidelity space in reduced-order space. Errors for HMC smoother are obtained for $100$ ensemble members with $25$ burn-in steps, and $5$ mixing steps.
The steps size for the symplectic integrator is empirically tuned and unified to $T=0.1$ with $h=0.01$, and $m=10$.
}
\includegraphics[width=10cm, height=6cm]%
{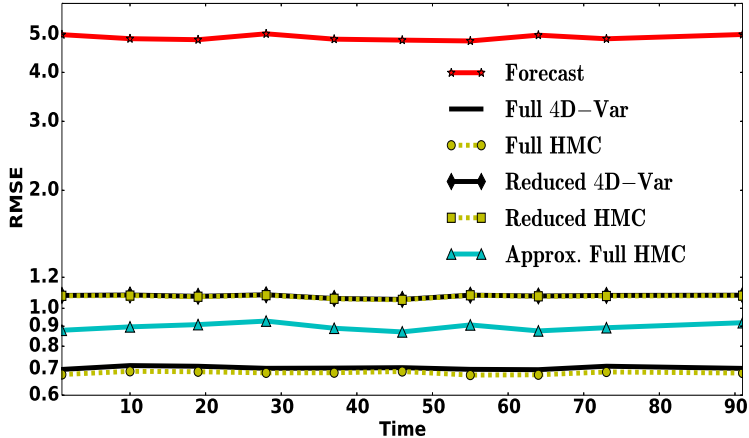}%
\label{fig:RMSE_Collective_Plot}
\end{figure}

The mean of the ensemble generated by HMC smoother is an MVUE of the posterior mean, and we are interested in comparing it against the 4D-Var solution. Figure~\ref{fig:RMSE_Collective_Plot} shows the Root mean squared (RMSE) errors associated with the 4D-Var and HMC estimates of the posterior mean. The size of the ensemble generated by the different HMC smoothers here is $\nens=100$. We see clearly that the MVUE generated by HMC in both full and reduced space is at least as good as the 4D-Var minimizer. It is obvious that using Algorithm~\ref{alg:HMC_Sampling_Smoother_combined} to sample the full space, while approximating the gradient using reduced-space information, results in an analysis that is better than the case where the sampler is limited to the reduced space.
In addition to testing the quality of the analysis (first-order moment here), we are interested in quantifying the quality of the analysis error covariance matrix generated by HMC. For reference we use HMC in full space to sample $\nens=1000$ members to produce a good estimate $ \A_0^{\rm ens} \approx \A_0$. In the cases of reduced space sampling and
approximate sampling in the full space we fix the ensemble size to $\nens=100$. To compare analysis error covariances obtained in the different scenarios we perform a statistical test of the hypothesis $\mathbf{H}_0:\Sigma_1=\Sigma_2$ for the equality of two covariance matrices. Since the state space dimension can be much larger than ensemble size, we choose the test statistic ~\cite{schott2007test}  that works in high dimensional settings. Assume we have two probability distributions with covariance matrices $\Sigma_1,\, \Sigma_2$ respectively, and consider sample estimates ${\Ss_1,\, \Ss_2}$ obtained using ensembles of sizes $n_1$ and $n_2$, respectively. The test statistic $t_{mn}^*$ defined in~\eqref{eqn:Scott_Test_Statistic} asymptotically follows a standard normal distribution in the limit of large ensemble size and state space dimension.
At a significance level $\alpha$, the two sided test $\mathbf{H}_0:\Sigma_1=\Sigma_2$ is rejected only if $| t_{mn}^* | > z_{\alpha/2}$, where $Z=\mathcal{N}(0,)$, and
$\PD(Z\ge z_{\alpha/2}) = \alpha/2$.
\begin{equation}
\label{eqn:Scott_Test_Statistic}
\begin{aligned}
t_{mn}^* &= \frac{t_{mn}}{\hat{\theta}}\,, \\
t_{mn}   &= {
 \left( 1-\frac{n_1-2}{\eta_1} \tr{\left(\Ss_1^2\right)} \right)
  +  \left( 1-\frac{n_2-2}{\eta_2} \tr{\left(\Ss_2^2\right)} \right)
  -2 \tr{\left(\Ss_1 \Ss_2\right)}
  - \frac{n_1}{\eta_1} \left( \tr{\left(\Ss_1\right)} \right)^2
  - \frac{n_2}{\eta_2} \left( \tr{\left(\Ss_2\right)} \right)^2
}\,,  \\
\eta_1 &= (n_1+2)(n_1-1),\quad \eta_2=(n_2+2)(n_2-1)\,,  \\
n &= n_1 + n_2,\quad \Ss = \frac{n_1}{n} \Ss_1 + \frac{n_2}{n} \Ss_2\,,  \\
\hat{\theta} &= \sqrt{4a^2 \left( \frac{n_1+n_2}{n_1 n_2} \right)^2 }\,, \quad
a = \frac{n^2}{(n+2)(n-1)} \left( \tr{(\Ss^2)} - \frac{(\tr{(\Ss)})^2}{n} \right)\,.
\end{aligned}
\end{equation}
Table~\ref{table:Equal_Covariances_Hypothesis_Tests} shows the results of the tests conducted to compare the covariance matrices.
\begin{table}[htbp]
\centering
\caption{Results of statistical tests conducted to compare covariance matrices obtained by HMC smoother in the three scenarios.
 $\Aini$ is the true posterior covariance of the distribution~\eqref{eqn:Smoothing_Posterior_Kernel}.
 $\widetilde{\A}_0$ is the true posterior covariance of the distribution with negative-log given by~\eqref{eqn:Reduced_HMC_Potential_Energy_and_Gradient}, while
 $\widetilde{\A}^{\rm ens}_0$ is the ensemble-based approximation obtained by Algorithm~\ref{alg:HMC_Sampling_Smoother_combined}.
 ${\wideparen \A}_0$ is the true posterior covariance of the distribution with negative-log given by~\eqref{eqn:Projected_HMC_Potential_Energy_and_Gradient}, while
 ${\wideparen \A}^{\rm ens}_0$ is the ensemble-based approximation obtained by Algorithm~\ref{alg:HMC_Sampling_Smoother_combined}.
 }
\label{table:Equal_Covariances_Hypothesis_Tests}
\begin{tabular}{|c|p{4cm}|c|c|c|c|} \hline
  \multicolumn{3}{|c|}{Test}  & Ensemble statistics  & Test-statistic \\ \hline
1 & Sampling the reduced space
  &\begin{tabular}[c]{@{}l@{}} $\mathbf{H_0}:\Aini=\widetilde{\A}_0$ \\ $\mathbf{H_a}:\Aini\neq\widetilde{\A}_0$ \end{tabular}
& \begin{tabular}[c]{@{}l@{}} $n_1 = 1000,\, n_2 = 100$ \\ $\Ss_1 = \Aini^{\rm ens},\, \Ss_2 = \widetilde{\A}_0^{\rm ens}$ \end{tabular}  & $t_{nm}^*=61.0258$  \\ \hline
2 & Sampling the full space with approximate gradient
  & \begin{tabular}[c]{@{}l@{}} $\mathbf{H_0}:\Aini={\wideparen{\A}}_0$ \\ $\mathbf{H_a}:\Aini\neq{\wideparen{\A}}_0$ \end{tabular}
& \begin{tabular}[c]{@{}l@{}} $n_1 = 1000,\, n_2 = 100$ \\ $\Ss_1 = \Aini^{\rm ens},\, \Ss_2 = {\wideparen \A}_0^{\rm ens}$ \end{tabular}   & $t_{nm}^*=2.4514$  \\ \hline
\end{tabular}
\end{table}
In the case of sampling in the reduced space the null hypothesis is rejected due to strong evidence based on the
samples' estimates. For the approximate full space sampling at a significance level $\alpha=0.01$ there is no significant evidence to supports rejection. This gives a strong indication that the ensemble generated in the second case describes the uncertainty in the analysis much better than the first case. The test results at least don't oppose the conclusion that sampling~\eqref{eqn:Full_PDF_with_Approximate_Gradient} using Algorithm~\ref{alg:HMC_Sampling_Smoother_combined} results in ensembles capable of  estimating
the posterior covariance matrix.
%
\subsection{Computational costs}
\begin{figure}
\centering
\caption{Data assimilation results using 4D-Var schems, and HMC smoother, in both high-fidelity space in reduced-order space.
 CPU-times for HMC smoother are obtained for $30$ ensemble members with $25$ burn-in steps, and $5$ mixing steps. The steps size
 for the symplectic integrator is empirically tuned and unified to $T=0.1$ with $h=0.01$, and $m=10$.
 The red color represents the CPU-time spent during optimization steps only. Blue and Green colors, respectively, represent CPU-time spent during the burn-in and
 the sampling( and mixing) steps. 
 }
\includegraphics[width=0.60\linewidth]%
    {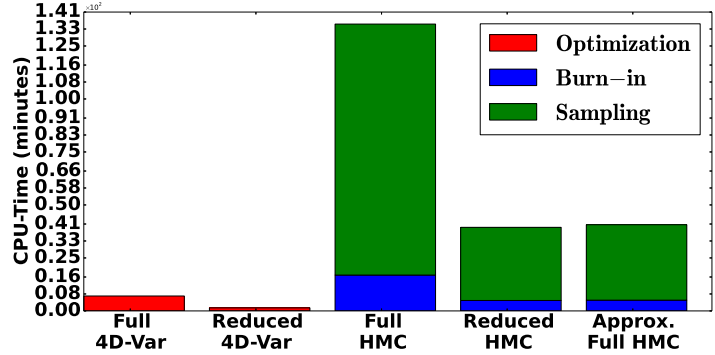}%
\label{fig:CPU_Times}
\end{figure}

%
The computational cost for HMC smoother in full space is much higher than the cost of 4D-Var, however it comes with the advantage of generating a consistent estimate of the analysis error covariance matrix. The bottleneck of HMC smoother is the propagation of the forward and backward model to evaluate the gradient of the potential energy. Using surrogate models radically reduces the computational cost. A detailed discussion of the computational cost in terms of model propagation
can be found in~\cite{Attia_HMCSmoother_TR}. Here we report the CPU time of the different scenarios as shown in Figure~\ref{fig:CPU_Times} and Table~\ref{table:CPU_Times}. The HMC CPU-time also depends on the settings of the parameters and the size of the ensemble. Following \cite{Attia_HMCSmoother_TR} we compare the CPU-times to generate $30$ ensemble members.
The CPU-times are almost similar when the two strategies in Algorithm~\ref{alg:HMC_Sampling_Smoother_combined} are applied, and both are approximately four times faster than the original HMC smoother. The online cost of the approximate smoother is still higher than the cost of 4D-Var in full space, however it is notably reduced by using information coming from a reduced space. The cost can be further reduced by cleverly tuning the sampler parameters or projecting the observation operator and observation error statistics in the reduced space. These ideas will be considered in the future to further reduce the cost of the HMC sampling smoother.
It is very important to highlight the fact that the goal is not just to find an anlysis state but to approximate the whole posterior distribution.
Despite the high cost of the HMC smoother, we obtain a consistent description of the uncertainty of the analysis state, e.g. an estimate of the posterior covariances.
\begin{table}[h]
\centering
\caption{Data assimilation results using 4D-Var schemes, and HMC smoother, in both high-fidelity space in reduced-order space.
 CPU-times for HMC smoother are obtained for $30$ ensemble members with $25$ burn-in steps, and $5$ mixing steps. The steps size
 for the symplectic integrator is empirically tuned and unified to $T=0.1$ with $h=0.01$, and $m=10$. }
\label{table:CPU_Times}
\resizebox{\columnwidth}{!}{%
\begin{tabular}{|p{1.6cm}|p{2cm}|p{2cm}|p{1.8cm}|p{1.1cm}|p{1.8cm}|p{1.1cm}|p{1.8cm}|p{1.1cm}|}
\hline
\multirow{4}{*}{\textbf{Cost}} & \multicolumn{8}{c|}{\textbf{Experiment}}\\ \cline{2-9}
   & \multicolumn{2}{c|}{\textbf{4D-Var}}  & \multicolumn{6}{c|}{\textbf{HMC Smoother}}  \\ \cline{2-9}
   & \multirow{2}{*}{\textbf{\begin{tabular}[c]{@{}l@{}} high-fidelity \\ space \end{tabular} }} & \multirow{2}{*}{\textbf{\begin{tabular}[c]{@{}l@{}} reduced-order\\ space \end{tabular} }}  & \multicolumn{2}{p{3cm}|}{\textbf{high-fidelity space} }& \multicolumn{2}{p{2.8cm}|}{\textbf{reduced-fidelity space} } & \multicolumn{2}{p{3cm}|}{\textbf{high-fidelity space with approximate gradient}} \\ \cline{4-9}
   &   & & \textit{average per ensemble member} & \textit{total} & \textit{average per ensemble member} & \textit{total} & \textit{average per ensemble member}& \textit{total}\\ \hline
\textbf{CPU-time (minutes)}  & $7.04$    & $1.44$  & $0.68$   & $118.42$   & $0.20$  & $34.50$     & $0.20$  & $35.58$ \\ \hline  
\end{tabular}
}
\end{table}
%

\section{Conclusions and Future Work}
\label{Sec:Conclusion}
%
The HMC sampling smoother is developed as a general ensemble-based data assimilation framework to solve the non-Gaussian four-dimensional data assimilation problem. The original formulation of the HMC smoother works with the full dimensional model. It provides a consistent description of the posterior distribution, however it is very expensive due to the necessary large number of full model runs. The HMC sampling smoother employs reduced-order approximations of the model dynamics. It achieves computational efficiency while retaining most of the accuracy of the full space HMC smoother. The formulations discussed here still assume a Gaussian prior at the initial time, which is a weak assumption since the forward propagation through nonlinear model dynamics will result in a non-Gaussian likelihood. This assumption, however, can be easily relaxed using a mixture of Gaussians to represent the background at the initial time; this will be considered in future work. We plan to explore the possibility of using the KL-Divergence measure between the high fideltity distribution
and both the projected and the approximate posterior distribution, to guide the optimal choice of the size of reduced-order basis. In future work we will also consider incorporating an HMC sampler capable of automatically tuning the parameters of the symplectic integrator, such as NUTS~\cite{hoffman2014NUTS}, in order to further enhance the smoother performance.

\section*{Acknowledgments}

The work of Dr. R\u azvan Stefanescu and Prof. Adrian Sandu was supported by awards NSF CCF--1218454, NSF DMS--1419003, AFOSR FA9550--12--1--0293--DEF, AFOSR 12-2640-06, and by the Computational Science Laboratory at Virginia Tech.

\bibliographystyle{plain} 

\bibliography{Bib/Software,Bib/ROM_state_of_the_art,Bib/comprehensive_bibliography1,Bib/Razvan_bib,Bib/Razvan_bib_ROM_IP,Bib/reduced_models,Bib/data_assim_HMC,Bib/data_assim_kalman,Bib/data_assim_fdvar}

\end{document}